\documentclass[final]{siamltex}

\usepackage{amsmath}
\usepackage{graphics}
\usepackage{enumerate}
\usepackage{epsfig}
\usepackage{amsfonts}
\usepackage{amssymb}
\usepackage{showkeys}

\newtheorem{remark}{Remark}

\def\d{\,\mathrm{d}}
\def\rho{\varrho}
\def\eps{\varepsilon}
\def\phi{\varphi}
\def\theta{\vartheta}
\def\N{\mathbb{N}}
\def\R{\mathbb{R}}
\def\C{\hbox{\rlap{\kern.24em\raise.1ex\hbox
      {\vrule height1.3ex width.9pt}}C}}
\def\P{\hbox{\rlap{I}\kern.16em P}}
\def\Q{\hbox{\rlap{\kern.24em\raise.1ex\hbox
      {\vrule height1.3ex width.9pt}}Q}}
\def\M{\hbox{\rlap{I}\kern.16em\rlap{I}M}}
\def\Z{\hbox{\rlap{Z}\kern.20em Z}}
\def\({\begin{eqnarray}}
\def\){\end{eqnarray}}
\def\[{\begin{eqnarray*}}
\def\]{\end{eqnarray*}}
\def\part#1#2{{\partial #1\over\partial #2}}

\def\grad{\nabla}
\def\Norm#1{\left\| #1 \right\|}
\def\pmb#1{\setbox0=\hbox{$#1$}
  \kern-.025em\copy0\kern-\wd0
  \kern-.05em\copy0\kern-\wd0
  \kern-.025em\raise.0433em\box0 }
\def\bar{\overline}

\def\tot#1#2{\frac{\d #1}{\d #2}} 

\def\d{\,\mathrm{d}}
\def\N{\mathbb{N}}
\def\R{\mathbb{R}}
\def\supp{{\rm supp\, }}
\def\rho{\varrho}
\def\eps{\varepsilon}
\def\phi{\varphi}

\def\O{\mathcal{O}}

\def\D{\mathcal{D}}
\def\E{\mathcal{E}}
\def\tx{\tilde x}

\title{Particle systems and kinetic equations modeling interacting agents in high dimension}

\author{M.~Fornasier\footnotemark[2]\ \footnotemark[3]
\and J.~Ha\v{s}kovec\footnotemark[3]\
\and J.~Vyb\'{i}ral\footnotemark[3]}

\begin{document}

\maketitle

\renewcommand{\thefootnote}{\fnsymbol{footnote}}

\footnotetext[2]{Faculty of Mathematics, Technical University of Munich, Boltzmannstrasse 3, D-85748 Garching, Germany}
\footnotetext[3]{Johann Radon Institute for Computational and
Applied Mathematics, Austrian Academy of Sciences, Altenbergerstrasse 69, A-4040 Linz, Austria}

\renewcommand{\thefootnote}{\arabic{footnote}}

\begin{abstract}
In this paper we explore how concepts of high-dimensional data compression via random projections onto lower-dimensional spaces 
can be applied for tractable simulation of certain dynamical systems modeling complex interactions. 
In such systems, one has to deal with a large number of agents (typically millions) 
in spaces of parameters describing each agent of high dimension (thousands or more).
Even with today's powerful computers, numerical simulations of such systems
are prohibitively expensive.
We propose an approach for the simulation
of dynamical systems governed by functions of {\it adjacency matrices} in high dimension, by random projections via Johnson-Lindenstrauss embeddings,
and recovery by compressed sensing techniques. We show how these
concepts can be generalized to work for associated kinetic equations, by addressing
the phenomenon of the delayed curse of dimension, known in information-based complexity for optimal numerical integration problems and measure quantization in high dimensions.
\end{abstract}



\pagestyle{myheadings}
\thispagestyle{plain}
\markboth{M.~FORNASIER, J.~HA\v{S}KOVEC AND J.~VYB\'{I}RAL}{PARTICLE AND KINETIC MODELING IN HIGH DIMENSION}

\section{Introduction}\label{Sec:Introduction}

The dimensionality scale of problems arising in our modern information
society has become very large and finding appropriate methods
for dealing with them is one of the great challenges of today's numerical simulation.
The most notable recent advances in data analysis are based
on the observation that in many situations, even for very complex
phenomena, the intrinsic dimensionality of the data is significantly lower than the ambient
dimension. Remarkable progresses have been made in data compression, processing, and acquisition.
We mention, for instance, the use of {\it diffusion maps} for data clouds and graphs in high dimension \cite{BN01,BN03,CLLMNWZ05a,CLLMNWZ05b,CL06,KMFK}
in order to define low-dimensional local representations of data with small distance distortion, and meaningful automatic clustering properties.
In this setting the embedding of data is performed by a {\it highly nonlinear} procedure, obtained by computing the eigenfunctions of suitable normalized
diffusion kernels, measuring the probability of transition from one data point to another over the graph.

Quasi-isometrical {\it linear} embeddings of high-dimensional point clouds into
low-dimen\-sional spaces of parameters are provided by the well-known Johnson-Lindenstrauss Lemma \cite{A,DG,JL}:
any cloud of $\mathcal N$ points in $\mathbb R^d$ can be embedded  by a random linear projection $M$ nearly
isometrically into $\mathbb R^k$ with $k = \mathcal O(\varepsilon^{-2} \log(\mathcal N))$
(a precise statement will be given below). This embedding strategy is simpler than the use of diffusion maps,
as it is linear, however it is ``blind'' to the specific geometry and local dimensionality of the data,
as the embedding dimension $k$ depends exclusively on the number of points in the cloud.
In many applications, this is sufficient, as the number of points $\mathcal N$ is supposed
to be a power of the dimension $d$, and the embedding produces an effective reduction
to $k = \mathcal O(\varepsilon^{-2} \log(\mathcal N)) = \mathcal O(\varepsilon^{-2} \log(d))$ dimensions.
As clarified in \cite{BDDW,KW}, the Johnson-Lindenstrauss Lemma is also at the basis of the possibility of performing optimal compressed and nonadaptive acquisition of high-dimensional data.
In {\it compressed sensing} \cite{carota06,do06-2} a vector $x \in \mathbb R^d$
is encoded in a vector $y\in\mathbb R^k$ by applying a random projection $M$,
which is modeling a linear acquisition device with random sensors, i.e., $y = Mx$.
From $y$ it is possible to decode $x$ approximately (see Theorem \ref{thmStrongRIP} below) by solving the convex optimization problem
$$
x^{\#}=\arg \min_{M z =y} \left ( \| z\|_{\ell_1^d} := \sum_{i=1}^d |z_i| \right),
$$
with the error distortion
$$
\| x^{\#} -x\|_{\ell_1^d} \leq C \sigma_K(x)_{\ell_1^d},
$$
where $\sigma_K(x)_{\ell_1^d} = \inf_{z: \#\supp(z) \leq K} \|z -x\|_{\ell_1^d}$ and $K = \mathcal O(k/(\log(d/k)+1))$. We denote $\Sigma_K=\{z \in \mathbb R^d: \#\supp(z) \leq K\}$ the set of $K$-sparse vectors, i.e., the union of $K$-dimensional coordinate subspaces in $\mathbb R^d$.
In particular, if $x \in \Sigma_K$, then $x^{\#} =x$.
Hence, not only is $M$ a Johnson-Lindenstrauss embedding, quasi-isometrical on point clouds and  $K$-dimensional coordinate subspaces, but also allows for the recovery of the most relevant components of high-dimensional vectors, from low-dimensional encoded information.
A recent work \cite{BW09,W08} extends the quasi-isometrical properties of the Johnson-Lindenstrauss embedding
from point clouds and  $K$-dimensional coordinate subspaces to smooth compact 
Riemannian manifolds with bounded curvature. Inspired by this work, in \cite{IM11} the authors extend the principles of compressed sensing in terms of point recovery on smooth compact Riemannian manifolds.

Besides these relevant results in compressing and coding-decoding high-dimensional ``stationary'' data, dimensionality reduction of complex dynamical systems and high-dimensional partial 
differential equations is a subject of recent intensive research. Several tools have been employed, for instance, the use of diffusion maps for dynamical systems \cite{NLK06}, 
tensor product bases and sparse grids for the numerical solution of linear high-dimensional PDEs \cite{DSS,BG,GK,GO}, the reduced basis method for solving high-dimensional 
parametric PDEs \cite{BCDDPW,BMPPTxx,MPT02,RHP08,S08,VPRP03}.
\\
In this paper we shall further explore the connection between data compression and {\it tractable} numerical simulation of dynamical systems. Eventually we address the solutions of associated high-dimensional kinetic equations.
We are specially interested in dynamical systems of the type 
\begin{equation}\label{eq:dyn}
\dot x_i(t)=f_i(\D x(t))+\sum_{j=1}^N f_{ij}(\D x(t))x_j(t),
\end{equation}
where we use the following notation:
\begin{itemize}
\item $N\in \N$ - number of agents,
\item $x(t)=(x_1(t),\dots,x_N(t))\in\R^{d\times N}$, where $x_i:[0,T]\to \R^d$, $i=1,\dots,N$,
\item $f_i:\R^{N\times N}\to \R^d,\quad i=1,\dots,N,$
\item $f_{ij}:\R^{N\times N}\to \R, \quad i,j=1,\dots,N$,
\item $\D:\R^{d\times N}\to \R^{N\times N}$, $\D x:=(\|x_i-x_j\|_{\ell_2^d})_{i,j=1}^N$ is the {\it adjacency matrix} of the point cloud $x$.
\end{itemize}
We shall assume that the governing functions $f_i$ and $f_{ij}$ are Lipschitz, but we shall specify the details later on.
The system \eqref{eq:dyn} describes the dynamics of multiple complex agents $x(t)=(x_1(t),\dots,x_N(t))\in\R^{d\times N}$, interacting on
the basis of their mutual ``social'' distance $\D x(t)$, and its general form includes several models for swarming and collective motion of animals and micro-organisms, aggregation of cells, etc. Several relevant effects can be included in the model by means of the functions $f_i$ and $f_{ij}$, in particular, fundamental binary mechanisms of
{\it attraction, repulsion, aggregation, self-drive}, and {\it alignment} \cite{CCR,CFTV,CS1,CS2,DOrsogna,KS}.
Moreover, possibly adding stochastic terms of random noise may also allow to consider {\it diffusion} effects \cite{BCC,CFTV}.
However, these models and motion mechanisms are mostly derived borrowing a leaf from physics, by assuming the agents (animals, micro-organisms, cells etc.)  as pointlike and exclusively determined by their spatial position and velocity in $\mathbb R^{d}$ for $d=3+3$.
In case we wished to extend such models of social interaction to more ``sophisticated'' agents,
described by many parameters ($d \gg 3 + 3$),
the simulation may become computationally prohibitive.
Our motivation for considering high-dimensional situations stems from the modern development of communication
technology and Internet, for which we witness the development of larger and larger communities accessing
information (interactive databases), services (financial market), social interactions (social networks) etc.
For instance, we might be interested to simulate the behavior of certain subsets of the financial market
where the agents are many investors, who are characterized by their portfolios of several hundreds of investments.
The behavior of each individual investor depends on the dynamics of others according to a suitable
social distance determined by similar investments.
Being able to produce meaningful simulations and learning processes of such complex dynamics is an issue,
which might be challenged by using suitable compression/dimensionality reduction techniques.\\
The idea we develop in this paper is to randomly project the system  and its initial condition by Johnson-Lindenstrauss embeddings
to a lower-dimensional space where an independent simulation can be performed with significantly reduced complexity.
We shall show that the use of multiple projections and parallel computations allows
for an approximate reconstruction of the high-dimensional dynamics, by means of compressed sensing techniques. After we explore the tractable simulation of the dynamical systems \eqref{eq:dyn} when the dimension $d$ of the parameter space is large, we also address the issue of whether we can perform tractable simulations when the number $N$ of agents is getting very large. Unlike the control of a finite number of agents, the numerical
simulation of a rather large population of interacting agents ($N \gg 0$) can constitute a serious difficulty which stems from the accurate solution of a possibly very large system of ODEs. Borrowing the
strategy from the kinetic theory of gases \cite{cip}, we want instead to consider a density distribution of agents, depending on their $d$-parameters, which interact with stochastic influence (corresponding to classical collisional rules
in kinetic theory of gases) -- in this case the influence is ``smeared'' since two individuals may interact also when they are far apart in terms of their ``social distance'' $\D x$. Hence, instead of simulating the behavior of
each individual agent, we shall describe the collective behavior encoded by a density distribution $\mu$, whose evolution is governed by one sole mesoscopic partial differential equation. We shall show that, under realistic assumptions on the concentration of the measure $\mu$ on sets of lower dimension, we can also acquire information on the properties of the high-dimensional measure solution $\mu$ of the corresponding kinetic equation, by considering random projections to lower dimension. Such approximation properties are determined by means of the combination of optimal numerical integration principles for the high-dimensional measure $\mu$ \cite{GL,Gruber} and the results previously achieved for particle dynamical systems.

\subsection{Fundamental assumptions}

We introduce the following notation for $\ell_p$-norms of vectors $v \in \mathbb R^d$,
$$
\| v\|_{\ell_p^d} := \left( \sum_{i=1}^d |v_i|^p \right )^{1/p}\qquad \mbox{for } 1\leq p < \infty,
$$
and 
$$
\| v\|_{\ell_\infty^d}:=\max_{i=1,\dots,d} |v_i|.
$$
For matrices $x \in \mathbb R^{n \times m}$ we consider the mixed norm
$$
\|x\|_{\ell_p^m(\ell_q^n)} := \| (\|x_i\|_{\ell_p^n})_{i=1}^m \|_{\ell_q^m},
$$
where $x_i \in \mathbb R^n$ is the $i^{th}$-column of the matrix $x$.

For the rest of the paper we impose three fundamental assumptions about
Lipschitz and boundedness properties of $f_i$ and $f_{ij}$,
\begin{align}
\label{eq:condf1} |f_i(a)-f_i(b)|&\le L\|a-b\|_{\ell_\infty^N(\ell_\infty^N)}, \quad i=1,\dots,N\\
\label{eq:condf2}\max_{i=1,\dots,N}\sum_{j=1}^N|f_{ij}(a)|&\le L',\\
\label{eq:condf3}\max_{i=1,\dots,N}\sum_{j=1}^N|f_{ij}(a)-f_{ij}(b)|&\le L''\|a-b\|_{\ell_\infty^N(\ell_\infty^N)},
\end{align}
for every $a,b\in\R^{N\times N}$.
Unfortunately, models of real-life phenomena would not always satisfy these conditions,
for instance models of financial markets or socio-economic interactions can be expected to exhibit severely discontinuous behavior.
However, these assumptions are reasonable in certain regimes and allow us to prove the concept
we are going to convey in this paper, i.e., the possibility of simulating high-dimensional dynamics
by multiple independent simulations in low dimension.

\subsection{Euler scheme, a classical result of stability and convergence, and its complexity}
We shall consider the system of ordinary differential equations of the form \eqref{eq:dyn} with the initial condition
\begin{equation}\label{eq:dyn2}
x_i(0) = x_i^0 \,,\qquad i=1,\dots,N \,.
\end{equation}
The Euler method for this system is given by \eqref{eq:dyn2} and 
\begin{equation}\label{eq:eul1}
x_i^{n+1}:=x_i^n+h\left[f_i(\D x^n)+\sum_{j=1}^N f_{ij}(\D x^n)x^n_j\right],\quad n=0,\dots,n_0-1.
\end{equation}
where $h>0$ is the time step and $n_0:=T/h$ is the number of iterations. We consider here the {\it explicit} Euler scheme exclusively for the sake of simplicity, for more sophisticated integration methods might be used.

The simulation of the dynamical system \eqref{eq:eul1} has a complexity which is at least the one of computing the {\it adjacency matrix} $\D \tilde x^n$ at each discrete  time $t^n$, i.e., $\mathcal O( d \times N^2)$.
The scope of the next sections is to show that, up to an $\varepsilon$-distortion, we can approximate the dynamics of \eqref{eq:dyn} by projecting the system into 
lower dimension and by executing in parallel computations with reduced complexity. Computation of the adjacency matrix in the new dimension requires only $\mathcal O( \varepsilon^{-2} \log (N)\times N^2)$ operations.
Especially if the distortion parameter $\varepsilon>0$ is not too small and the number of agents is of a polynomial order in $d$, we reduce the complexity of computing the adjacency matrix to $\mathcal O(\log (d)\times N^2)$.

\section{Projecting the Euler method: dimensionality reduction of discrete dynamical systems}\label{Sec:Euler}

\subsection{Johnson-Lindenstrauss embedding}

We wish to project the dynamics of \eqref{eq:dyn} into a lower-dimensional space
by employing a well-known result of Johnson and Lindenstrauss \cite{JL},
which we informally rephrase for our purposes as follows.

\begin{lemma}[Johnson and Lindenstrauss]
\label{JLlem}
Let $\mathcal P$ be an arbitrary set of $\mathcal N$ points in $\R^d$.
Given a distortion parameter $\varepsilon>0$, there exists a constant
$$k_0= \mathcal O(\varepsilon^{-2}\log (\mathcal N)),$$
such that for all integers $k\ge k_0$, there exists a $k\times d$ matrix $M$ for which 
\begin{equation}\label{eq:JLlem}
(1-\varepsilon)\|x-\tx\|_{\ell_2^d}^2\le \|Mx - M\tx\|_{\ell_2^k}^2\le (1+\varepsilon)\|x-\tx\|_{\ell_2^d}^2,
\end{equation}
 for all $x,\tx\in \mathcal P$.
\end{lemma}
It is easy to see that the condition
\begin{equation}\label{eq:RIP1}
(1-\varepsilon)\|p\|_{\ell_2^d}^2\le \|M p\|_{\ell_2^k}^2\le (1+\varepsilon)\|p\|_{\ell_2^d}^2, \quad p \in \mathbb R^d,
\end{equation}
implies
\begin{equation}\label{eq:RIP2}
(1-\varepsilon)\|p\|_{\ell_2^d} \le \|M p\|_{\ell_2^k} \le (1+\varepsilon)\|p\|_{\ell_2^d}, \quad p \in \mathbb R^d,
\end{equation}
for $0<\varepsilon<1$, which will be used in the following sections. On the other hand, \eqref{eq:RIP2}
implies \eqref{eq:RIP1} with $3\varepsilon$ instead of $\varepsilon$. 

Our aim is to apply this lemma to dynamical systems. As the mapping $M$ from Lemma \ref{JLlem}
is linear and almost preserves distances between the points (up to the $\varepsilon>0$ distortion as described above),
we restrict ourselves to dynamical systems which are quasi-linear or whose non-linearity depends only
on the mutual distances of the points involved, as in \eqref{eq:dyn}.

Let us define the additional notation, which is going to be fixed throughout the paper:

\begin{itemize}
\item $d\in \N$ - dimension (large),
\item $\varepsilon>0$ - the distortion parameter from Lemma \ref{JLlem},
\item $k\in \N$ - new dimension (small),
\item $M\in \R^{k\times d}$ - randomly generated matrix as described below.
\end{itemize}

The only constructions of a matrix $M$ as in Lemma \ref{JLlem} known up to now are stochastic, i.e., the matrix is 
randomly generated and has the 
quasi-isometry property \eqref{eq:JLlem} with high probability.
We refer the reader to \cite{DG} and \cite[Theorem 1.1]{A} for two typical versions of the Johnson-Lindenstrauss Lemma.

We briefly collect below some well-known instances of random matrices, which satisfy the statement of Lemma \ref{JLlem} with high probability:
\begin{itemize}
\item $k \times d$ matrices $M$ whose entries $m_{i,j}$ are independent realizations of Gaussian random variables
$$
m_{i,j} \sim \mathcal N\left (0, \frac{1}{k} \right );
$$
\item $k \times d$ matrices $M$ whose entries are independent realizations of $\pm$ Bernoulli random variables
$$
m_{i,j}:=\left \{
\begin{array}{ll}
+\frac{1}{\sqrt k},& \text{ with probability } \frac{1}{2}\\
-\frac{1}{\sqrt k},& \text{ with probability } \frac{1}{2}
\end{array}
\right .
$$
\end{itemize}

Several other random projections suitable for Johnson-Lindenstrauss embeddings can be constructed following Theorem \ref{KWthm} recalled below, and we refer the reader to \cite{KW} for more details. 

\subsection{Uniform estimate for a general model}

If $M\in \R^{k\times d}$ is a matrix, we consider the projected Euler method in $\R^k$ associated to the high-dimensional system \eqref{eq:dyn2}-\eqref{eq:eul1}, namely
\begin{align}
\label{eq:eul2}y_i^0&:=Mx_i^0,\\
\label{eq:eul3}y_i^{n+1}&:=y_i^n+h\left[Mf_i(\D' y^n)+\sum_{j=1}^N f_{ij}(\D' y^n)y^n_j\right],\quad n=0,\dots,n_0-1.
\end{align}
We denote here  $\D':\R^{k\times N}\to \R^{N\times N}$, $\D' y:=(\|y_i-y_j\|_{\ell_2^k})_{i,j=1}^N$, the {\it adjacency matrix} of the agents $y=(y_1,\dots,y_N)$ in $\mathbb R^{k \times N}$.
The first result of this paper reads as follows.
\begin{theorem}\label{thm1}
Let the sequences 
$$
\{x_i^n, i=1,\dots,N \ \text{and}\ n=0,\dots,n_0\} \quad \text{and}\quad \{y_i^n, i=1,\dots,N \ \text{and}\ n=0,\dots,n_0\}
$$
be defined by \eqref{eq:dyn2}-\eqref{eq:eul1} and \eqref{eq:eul2}-\eqref{eq:eul3} with $f_i$ and $f_{ij}$ satisfying \eqref{eq:condf1}--\eqref{eq:condf3}
and a matrix $M\in\R^{k\times d}$ with
\begin{gather}
\label{eq:condM1}\left\|Mf_i(\D' y^n)-M f_i(\D x^n)\right\|_{\ell_2^k}\le (1+\varepsilon) \left\|f_i(\D' y^n)-f_i(\D x^n)\right\|_{\ell_2^d},\\
\label{eq:condM2}\|Mx_j^n\|_{\ell_2^k}\le (1+\varepsilon) \|x_j^n\|_{\ell_2^d},\\
\label{eq:condM3}(1-\varepsilon) \|x_i^n-x_j^n\|_{\ell_2^d}\le \|Mx_i^n-Mx_j^n\|_{\ell_2^k}\le (1+\varepsilon) \|x_i^n-x_j^n\|_{\ell_2^d}
\end{gather}
for all $i,j=1,\dots, N$ and all $n=0,\dots,n_0$. Moreover, let us assume that
\begin{align}
\alpha&\ge \max_j\|x^n_j\|_{\ell_2^d}\quad \text{for all}\quad n=0,\dots,n_0,\quad j=1,\dots,N. \label{alpha1}
\end{align}
Let 
\begin{equation}\label{eq:defe}
e_i^n:=\|y_i^n-Mx_i^n\|_{\ell_2^k},\ i=1,\dots,N \ \text{and}\ n=0,\dots,n_0
\end{equation}
and set $\E^n:=\max_ie_i^n$. Then
\begin{equation}\label{eq:unifstab}
{\mathcal E}^n\le \varepsilon hn B \exp(hnA),
\end{equation}
where $A:=L'+2(1+\varepsilon)(L+\alpha L'')$ and $B:=2\alpha(1+\varepsilon)(L+\alpha L'')$.
\end{theorem}

We remark that conditions \eqref{eq:condM1}-\eqref{eq:condM3} are in fact satisfied
as soon as $M$ is a suitable Johnson-Lindenstrauss embedding as in Lemma \ref{JLlem},
for the choice $\mathcal N = 2 N n_0$ and $k=\mathcal O(\varepsilon^{-2} \log(\mathcal N))$.

\begin{proof}
Using \eqref{eq:defe} and \eqref{eq:dyn2}-\eqref{eq:eul1} and \eqref{eq:eul2}-\eqref{eq:eul3} combined with \eqref{eq:condM1} and \eqref{eq:condM2},
we obtain
\begin{align*}
e_i^{n+1}&\le e_i^n+h\left \|Mf_i(\D' y^n)-M f_i(\D x^n)\right\|_{\ell_2^k}
+h\left\|\sum_{j=1}^N f_{ij}(\D' y^n)y_j^n-f_{ij}(\D x^n)Mx_j^n\right\|_{\ell_2^k}\\
&\le e_i^{n}+h(1+\varepsilon)\left\|f_i(\D' y^n)-f_i(\D x^n)\right\|_{\ell_2^d}\\
&\quad +h\sum_{j=1}^N \Bigl(\|f_{ij}(\D' y^n)y_j^n-f_{ij}(\D' y^n)Mx_j^n\|_{\ell_2^k}+\|f_{ij}(\D' y^n)Mx_j^n-f_{ij}(\D x^n)Mx_j^n\|_{\ell_2^k}\Bigr)\\
&\le e_i^{n}+h(1+\varepsilon)\left\|f_i(\D' y^n)-f_i(\D x^n)\right\|_{\ell_2^d}\\
&\quad +h\sum_{j=1}^N \Bigl(|f_{ij}(\D' y^n)|e_j^n+(1+\varepsilon)\|x_j^n\|_{\ell_2^d}\cdot|f_{ij}(\D' y^n)-f_{ij}(\D x^n)|\Bigr).
\end{align*}

Taking the maximum on both sides, this becomes
\begin{align*}
\E^{n+1}&\le \E^n +h(1+\varepsilon) \max_i \|f_i(\D' y^n)-f_i(\D x^n)\|_{\ell_2^d}\\
&\quad +h\E^n\max_i\sum_{j=1}^N |f_{ij}(\D' y^n)|+h(1+\varepsilon)\alpha \cdot \max_i \sum_{j=1}^N|f_{ij}(\D' y^n)-f_{ij}(\D x^n)|.
\end{align*}
We use \eqref{eq:condf1}--\eqref{eq:condf3} for $a=\D'y^n$ and $b=\D x^n$ to estimate all the terms on the right-hand side.
This gives
\begin{align*}
\E^{n+1}&\le \E^n +h(1+\varepsilon) L \|\D' y^n- \D x^n\|_{\ell_\infty^N(\ell_\infty^N)}+h\E^nL'+h(1+\varepsilon)\alpha L''\|\D' y^n- \D x^n\|_{\ell_\infty^N(\ell_\infty^N)}\\
&\le \E^n(1+hL')+h(1+\varepsilon)(L+\alpha L'')\left[\|\D' y^n-\D' Mx^n\|_{\ell_\infty^N(\ell_\infty^N)} + \|\D' M x^n-\D x^n\|_{\ell_\infty^N(\ell_\infty^N)}\right]\\
&\le \E^n(1+hL')+2h(1+\varepsilon)(L+\alpha L'')(\E^n+\alpha\varepsilon),
\end{align*}
where we used \eqref{eq:condM3} in the last line. This, together with ${\mathcal E}^0=0$, leads to
$$
{\mathcal E}^n\le \varepsilon hn B \exp(hnA),
$$
where $A:=L'+2(1+\varepsilon)(L+\alpha L'')$ and $B:=2\alpha(1+\varepsilon)(L+\alpha L'')$.
\end{proof}
\subsection{Uniform estimate for the Cucker-Smale model}\label{Subs:UniformCS}
As a relevant example, let us now show that Theorem \ref{thm1} can be applied to the well-known Cucker-Smale model, introduced and analyzed in \cite{CS1,CS2}, which is described by
\begin{align}
\dot x_i&=v_i \in \mathbb R^d, \label{CS1} \\
\dot v_i&=\frac{1}{N}\sum_{j=1}^N g(\|x_i-x_j\|_{\ell_2^d})(v_j-v_i), \quad i=1,\dots,N. \label{CS2}
\end{align}
The function $g:[0,\infty)\to \R$ is given by $g(s)=\frac{G}{(1+s^2)^{\beta}}$, for $\beta>0$, and bounded by $g(0)=G>0.$
This model describes the {\it emerging of consensus} in a group of interacting agents, trying to
{\it align} (also in terms of abstract consensus) with their neighbors. One of the motivations of the model 
from Cucker and Smale was to describe the formation and evolution of languages \cite[Section 6]{CS2}, although, due to its simplicity,
it has been eventually related mainly to the description of the {\it emergence of flocking} in groups of birds \cite{CS1}. In the latter case,
in fact, spatial and velocity coordinates are sufficient to describe a pointlike agent ($d=3+3$), while for the evolution of languages,
one would have to take into account a much broader dictionary of parameters,
hence a higher dimension $d\gg 3+3$ of parameters, which is in fact the case of our interest in the present paper.

Let us show that the model is indeed of the type \eqref{eq:dyn}. We interprete the system as a group of $2N$ agents in $\R^d$,
whose dynamics is given by the following equations
\begin{align*}
\dot x_i&=\sum_{j=1}^N f^{x}_{ij}v_j \in \mathbb R^d, \\
\dot v_i&=\sum_{j=1}^N f^{v}_{ij}(\D x)v_j, \quad i=1,\dots,N
\end{align*}
with $f^{x}_{ij}:=\delta_{ij}$, $\displaystyle f^v_{ii}(\D x):=-\frac{1}{N}\sum_{k=1}^N g(\|x_{i}-x_{k}\|_{\ell_2^d})$,
and $\displaystyle f^v_{ij}(\D x):=\frac{1}{N} g(\|x_{i}-x_{j}\|_{\ell_2^d})$, for $i\neq j$.
The condition \eqref{eq:condf1} is empty, \eqref{eq:condf2} reads
$$
L'\ge\max(1,2G)\ge \max_i\left\{1,\frac{2}{N}\sum_{k=1}^Ng(\|x^n_{i}-x^n_{k}\|_{\ell_2^d})\right\}.
$$
Finally, 
\begin{align*}
\max_i\frac{2}{N}\sum_{j=1}^N &\Bigl|g(\|x^n_{i}-x^n_{j}\|_{\ell_2^d})-g(\|y^n_{i}-y^n_{j}\|_{\ell_2^k})\Bigr|\\
&\le \max_i\frac{2\|g\|_{{\rm Lip}}}{N}\cdot \sum_{j=1}^N \Bigl| \|x^n_{i}-x^n_{j}\|_{\ell_2^d}-\|y^n_{i}-y^n_{j}\|_{\ell_2^k}\Bigr|\\
&\le 2\|g\|_{{\rm Lip}}\cdot \|\D' y^n-\D x^n\|_{\ell_\infty^N(\ell_\infty^N)}
\end{align*}
shows that $L''\le 2\|g\|_{{\rm Lip}}.$
The boundedness of the trajectories in the phase-space of \eqref{CS1}-\eqref{CS2}
at finite time has been proved, for instance, in \cite{HL}, see also \cite[Theorem 4.6]{CCR}.
The boundedness at finite time is clearly sufficient to define the constant $\alpha$
appearing in Theorem \ref{thm1}, also because we are mainly interested in the dynamics for short time, due
to the error propagation. Of course the constant $\alpha$ might grow with time, but, for instance,
for the Cucker-Smale system it grows at most linearly in time \cite{CFTV};
as in the error estimate \eqref{eq:unifstab}
we have an exponential function in time appearing, the possible linear growth can be considered a negligible issue; 
moreover, as our numerical experiments show, see Section \ref{Sec:Numerics}, the situation is much better in practice, 
and suitable scaling, as indicated below, allows us to assume in several circumstances that
the constant $\alpha$ is uniformly bounded for all times.
In fact, even when we were interested in longer time or even asymptotical behavior, especially when pattern 
formation is expected, then we would observe the following additional facts:
In the Cucker-Smale model the center of mass and the mean velocity are invariants of the dynamics.
Moreover the rate of communication between particles is given by $g(s)=\frac{G}{(1+s^2)^{\beta}}$.
When $\beta \leq 1/2$ it is know (see \cite{CFTV}) that the dynamics will converge to a flocking configuration.
In this case one can translate at the very beginning the center of mass and the mean velocity to $0$,
and the system will keep bounded for all times. Hence in this case the constant $\alpha$ can 
also be considered uniform for all times (not only bounded at finite time).

\subsection{Least-squares estimate of the error for the Cucker-Smale model}
The formula \eqref{eq:unifstab} provides the estimate of the maximum of the individual errors, 
i.e., $\mathcal E^n:=\| (y^n_i - M x^n_i)_{i=1}^N \|_{\ell_\infty^N(\ell_2^k)}$.
In this section we address the stronger $\ell_2^N(\ell_2^k)$-estimate for the error.
For generic dynamical systems \eqref{eq:dyn} such estimate is not available in general,
and one has to perform a case-by-case analysis. 
As a typical example of how to proceed, we restrict ourselves to the Cucker-Smale model,
just recalled in the previous section.
The forward Euler discretization of~\eqref{CS1}--\eqref{CS2} is given by
\begin{align}
\label{eq:l2:1} x^{n+1}_i&=x^n_i+h v^n_i,\\
\notag v^{n+1}_i&=v^n_i+\frac{h}{N}\sum_{j=1}^N g(\|x^n_i-x^n_j\|_{\ell_2^d})(v^n_j-v^n_i)
\end{align}
with initial data $x_i^0$ and $v_i^0$ given.
Let $M$ be again a suitable random matrix in the sense of Lemma \ref{JLlem}.
The Euler method of the projected system is given by the initial conditions
$y^0_i=Mx^0_i$ and $w^0_i=Mv^0_i$ and the formulas
\begin{align}
\label{eq:l2:2}y^{n+1}_i&=y^n_i+hw^n_i,\\
\notag w^{n+1}_i&=w^n_i+\frac{h}{N}\sum_{j=1}^N g(\|y^n_i-y^n_j\|_{\ell_2^k})(w^n_j-w^n_i).
\end{align}

We are interested in the estimates of the following quantities
\begin{align}
\label{eq:l2:e1}e_{x,i}^n&:=\|y_i^n-Mx_i^n\|_{\ell_2^k},\quad \E_{x}^n:=\sqrt{\frac{1}{N}\sum_{i=1}^N (e_{x,i}^n)^2}=
\frac{\|(y_i^n-Mx_i^n)_{i=1}^N\|_{\ell_2^N(\ell_2^k)}}{\sqrt N},\\
\label{eq:l2:e2}e_{v,i}^n&:=\|w_i^n-Mv_i^n\|_{\ell_2^k},\quad \E_{v}^n:=\sqrt{\frac{1}{N}\sum_{i=1}^N (e_{v,i}^n)^2}=
\frac{\|(w_i^n-Mv_i^n)_{i=1}^N\|_{\ell_2^N(\ell_2^k)}}{\sqrt N}.
\end{align}

\begin{theorem}
Let the sequences $\{x_i^n\}, \{v_i^n\}, \{y_i^n\}$, $\{w_i^n\}$, $\{e^n_{x,i}\}$ and $\{e^n_{v,i}\}$, $i=1,\dots,N$ and $n=1,\dots,n_0$
be given by \eqref{eq:l2:1}, \eqref{eq:l2:2}, \eqref{eq:l2:e1} and \eqref{eq:l2:e2}, respectively.
Let $\varepsilon>0$ and let us assume, that the matrix $M$ satisfies
\begin{align*}
(1-\varepsilon)\|x_i^n-x_j^n\|_{\ell_2^d}\le \|Mx_i^n-Mx_j^n\|_{\ell_2^k} &\le (1+\varepsilon) \|x_i^n-x_j^n\|_{\ell_2^d}\quad\text{and}\\
(1-\varepsilon) \|v_i^n-v_j^n\|_{\ell_2^d}\le \|Mv_i^n-Mv_j^n\|_{\ell_2^k} &\le (1+\varepsilon) \|v_i^n-v_j^n\|_{\ell_2^d}
\end{align*}
for all $i,j=1,\dots,N$ and $n=0,\dots,n_0$.

Then the error quantities $\E_{x}^n$ and $\E_{y}^n$ introduced in \eqref{eq:l2:e1} and \eqref{eq:l2:e2} satisfy
\begin{equation}\label{eq:leastsqstab}
\sqrt{({\mathcal E_x^{n}})^2+({\mathcal E_v^{n}})^2}
\le \varepsilon(1+\varepsilon)hn\|g\|_{\rm Lip}VX\exp(hn\|{\mathcal A}\|),
\end{equation}
where $V:=\max_{i,j,n}\|v^n_i-v^n_j\|_{\ell_2^d}$, $X:=\max_{i,j,n}\|x^n_i-x^n_j\|_{\ell_2^d}$ and
$$
{\mathcal A}=\left(\begin{matrix}0&1\\2(1+\varepsilon)\|g\|_{\rm Lip}V&2G\end{matrix}\right).
$$

\end{theorem}

\begin{proof}
Using \eqref{eq:l2:1} and \eqref{eq:l2:2}, we obtain
$$
e_{x,i}^{n+1}\le e_{x,i}^n+he_{v,i}^n\quad \text{and}\quad \E_x^{n+1}\le \E_x^n+h\E_v^n.
$$
To bound the quantity $\E^{n}_v$ we have to work more.
We add and subtract the term $g(\|y_i^n-y_j^n\|_{\ell_2^k})(Mv_j^n-Mv_i^n)$ and apply \eqref{eq:l2:1} and \eqref{eq:l2:2}. This leads to
\begin{align}
\notag e_{v,i}^{n+1}&\le e_{v,i}^n+\frac{h}{N}\sum_{j=1}^N \Bigl(\|g(\|y_i^n-y_j^n\|_{\ell_2^k})(w_j^n-w_i^n)\pm g(\|y_i^n-y_j^n\|_{\ell_2^k})(Mv_j^n-Mv_i^n)\\
\notag &\qquad-g(\|x_i^n-x_j^n\|_{\ell_2^d})(Mv_j^n-Mv_i^n)\|_{\ell_2^k}\Bigr)\\
\label{eq:l2:3}&\le e_{v,i}^n+\frac{h}{N}\sum_{j=1}^N g(\|y_i^n-y_j^n\|_{\ell_2^k})(e_{v,j}^n+e_{v,i}^n)\\
\notag &\qquad+\frac{(1+\varepsilon)h\|g\|_{\rm Lip}}{N}\cdot\sum_{j=1}^N\|v_j^n-v_i^n\|_{\ell_2^d}\cdot\bigl | \|x_i^n-x_j^n\|_{\ell_2^d}-\|y_i^n-y_j^n\|_{\ell_2^k}\bigr|.
\end{align}
We estimate the first summand in \eqref{eq:l2:3} 
$$
\frac{h}{N}\sum_{j=1}^N g(\|y_i^n-y_j^n\|_{\ell_2^k})(e_{v,j}^n+e_{v,i}^n)\le \frac{hG}{N}\bigl[Ne_{v,i}^n+\sum_{j=1}^Ne_{v,j}^n\bigr]
=hG e_{v,i}^n+\frac{hG}{N}\sum_{j=1}^Ne_{v,j}^n
$$
and its $\ell_2$-norm with respect to $i$ by H\"older's inequality
\begin{equation}\label{eq:l2:4}
h\sqrt{N}G\E_v^n+\frac{hG}{N}\Biggl(\sum_{i=1}^N\biggl(\sum_{j=1}^Ne_{v,j}^n\biggr)^2\Biggr)^{1/2} \le 2h\sqrt{N}G\E_v^n.
\end{equation}

To estimate the second summand in \eqref{eq:l2:3} we make use of
\begin{align*}
\bigl| \|x_i^n&-x_j^n\|_{\ell_2^d}-\|y_i^n-y_j^n\|_{\ell_2^k}\bigr|\\
&\le \bigl| \|x_i^n-x_j^n\|_{\ell_2^d}-\|Mx_i^n-Mx_j^n\|_{\ell_2^k}\bigr| + \bigl |\|Mx_i^n-Mx_j^n\|_{\ell_2^k}-\|y_i^n-y_j^n\|_{\ell_2^k}\bigr|\\
&\le \varepsilon \|x_i^n-x_j^n\|_{\ell_2^d}+e_{x,i}^n+e_{x,j}^n.
\end{align*}
We arrive at
\begin{align*}
&\frac{ (1+\varepsilon)h \|g\|_{\rm Lip}}{N}\sum_{j=1}^N \|v_j^n-v_i^n\|_{\ell_2^d}(\varepsilon \|x_i^n-x_j^n\|_{\ell_2^d}+e_{x,i}^n+e_{x,j}^n)\\
&\qquad\le \frac{(1+\varepsilon)h \|g\|_{\rm Lip}V}{N}
\biggl\{\varepsilon \sum_{j=1}^N\|x_i^n-x_j^n\|_{\ell_2^d}
+ N e_{x,i}^n + \sum_{j=1}^N e_{x,j}^n\biggr\}.
\end{align*}
The $\ell_2$-norm of this expression with respect to $i$ is bounded by
\begin{align}
\notag &\frac{(1+\varepsilon)h \|g\|_{\rm Lip}V}{N}\left\{\varepsilon\Bigl(\sum_{i=1}^N \Bigl(\sum_{j=1}^N\|x_i^n-x_j^n\|_{\ell_2^d}\Bigr)^2\Bigr)^{1/2}
+N\Bigl(\sum_{i=1}^N (e_{x,i}^n)^2\Bigr)^{1/2} + \sqrt{N}\sum_{j=1}^N e_{x,j}^n\right\}\\
\label{eq:l2:5}&\le (1+\varepsilon)h \|g\|_{\rm Lip}V\sqrt N(\varepsilon X+2\E_x^n).
\end{align}
Combining \eqref{eq:l2:3} with \eqref{eq:l2:4} and \eqref{eq:l2:5} leads to the recursive estimate
\begin{align}
\label{eq:finest}\E_x^{n+1}&\le \E_x^n+h\E_v^n,\\
\notag\E_v^{n+1}&\le \E_v^n+2hG\E_v^n+ h (1+\varepsilon) \|g\|_{\rm Lip}V\left\{\varepsilon X+2\E_x^n\right\},
\end{align}
which we put into the matrix form
\begin{equation}\label{eq:l2:6}
\left(\begin{matrix}\E_x^{n+1}\\ \E_v^{n+1}\end{matrix}\right)
\le{\mathcal A'}\left(\begin{matrix}\E_x^{n}\\\E_v^{n}\end{matrix}\right)
+\left(\begin{matrix}0\\(1+\varepsilon)\varepsilon h\|g\|_{\rm Lip}V X\end{matrix}\right),
\end{equation}
where ${\mathcal A'}$ is a $2\times 2$ matrix given by 
$$
{\mathcal A'}={\mathcal Id}+h{\mathcal A}:=\left(\begin{matrix}1&0\\0&1\end{matrix}\right)
+h\left(\begin{matrix}0&1\\2(1+\varepsilon)\|g\|_{\rm Lip}V&2G\end{matrix}\right).
$$

Taking the norms on both sides of \eqref{eq:l2:6} leads to
$$
\sqrt{({\mathcal E_x^{n+1}})^2+({\mathcal E_v^{n+1}})^2}\le (1+h\|{\mathcal A}\|) 
\sqrt{({\mathcal E_x^{n}})^2+({\mathcal E_v^{n}})^2}+\varepsilon(1+\varepsilon)h\|g\|_{\rm Lip}VX,
$$
which gives the least-squares error estimate \eqref{eq:leastsqstab}.
\end{proof}

\section{Dimensionality reduction for continuous dynamical systems}\label{Sec:Cont}

\subsection{Uniform estimates for continuous dynamical systems}

In this section we shall establish the analogue of the above results
for the continuous time setting of dynamical systems of the type~\eqref{eq:dyn},
\(    \label{Cont1}
   \dot x_i &=& f_i(\D x) + \sum_{j=1}^N f_{ij}(\D x)x_j \,,\qquad i=1,\dots,N \,,\\
     x_i(0) &=& x_i^0 \,,\qquad i=1,\dots,N \,.  \label{Cont2}
\)
We adopt again the assumptions about Lipschitz continuity and boundedness of the right-hand side made in Section~\ref{Sec:Euler},
namely~\eqref{eq:condf1}, \eqref{eq:condf2}, and \eqref{eq:condf3}.

\begin{theorem}\label{Thm:Cont}
Let $x(t)\in\R^{d\times N}$, $t\in[0,T]$, be the solution of the system~\eqref{Cont1}--\eqref{Cont2}
with $f_i$'s and $f_{ij}$'s satisfying~\eqref{eq:condf1}--\eqref{eq:condf3},
such that
\( \label{alpha}
   \max_{t\in[0,T]} \max_{i,j} \|x_i(t) - x_j(t)\|_{\ell_2^d} \leq \alpha \,.
\)
Let us fix $k\in\mathbb{N}$, $k \leq d$, and a matrix $M\in\R^{k\times d}$ such that
\(    \label{Cont-JL}
    (1-\eps) \|x_i(t) - x_j(t)\|_{\ell_2^d} \leq \|Mx_i(t) - Mx_j(t)\|_{\ell_2^k} \leq (1+\eps) \|x_i(t) - x_j(t)\|_{\ell_2^d} \,,    
\)
for all $t\in[0,T]$ and $i$, $j=1,\dots,N$.
Let $y(t)\in\R^{k\times N}$, $t\in[0,T]$ be the solution of the projected system
\begin{eqnarray}
    \dot y_i &=& Mf_i(\D' y) + \sum_{j=1}^N f_{ij}(\D' y) y_j \,,\qquad i=1,\dots, N\,, \nonumber\\
      y_i(0) &=& Mx_i^0 \,,\qquad i=1,\dots, N\,, \label{lowdimsys}
\end{eqnarray}
such that for a suitable $\beta>0$,
\(   \label{beta}
    \max_{t\in[0,T]} \Norm{y(t)}_{\ell_\infty^N(\ell_2^d)} \leq \beta  \,.
\)

Let us define the column-wise $\ell_2$-error $e_i(t) := \|y_i-Mx_i\|_{\ell_2^k}$ for $i=1,\dots,N$
and
\[
    \E(t) := \max_{i=1,\dots,N} e_i(t) = \Norm{y-Mx}_{\ell_\infty^N(\ell_2^k)} \,.
\]
Then we have the estimate
\begin{equation}
\label{eq:unifstab2}
    \E(t) \leq \eps\alpha t (L\Norm{M} + L''\beta)
          \exp\left[(2L\Norm{M} + 2\beta L'' + L') t\right]  \,.
\end{equation}
\end{theorem}

\begin{proof}
Due to~\eqref{eq:condf1}--\eqref{eq:condf3}, we have for every $i=1,\dots, N$ the estimate
\[
\tot{}{t} e_i &=& \frac{\langle y_i-Mx_i,\tot{}{t}(y_i-Mx_i)\rangle}{\|y_i-Mx_i\|_{\ell_2^k}}\le \left\|\tot{}{t}(y_i-Mx_i)\right\|_{\ell_2^k}\\
       &\leq&  \|M f_i(\D' y) - M f_i(\D x)\|_{\ell_2^k} + \sum_{j=1}^N \|f_{ij}(\D' y)y_j - f_{ij}(\D x)Mx_j\|_{\ell_2^k} \\
       &\leq&  L \Norm{M} \Norm{\D' y - \D x}_{\ell_\infty^N(\ell_\infty^N)}
            + \sum_{j=1}^N \left( \|f_{ij}(\D x)(Mx_j - y_j)\|_{\ell_2^k} + \|(f_{ij}(\D x)-f_{ij}(\D' y)) y_j\|_{\ell_2^k} \right) \\
       &\leq&   L \Norm{M} \Norm{\D' y - \D x}_{\ell_\infty^N(\ell_\infty^N)}
            + L' \Norm{Mx - y}_{\ell_\infty^N(\ell_2^k)} + L'' \Norm{\D x-\D' y}_{\ell_\infty^N(\ell_\infty^N)} \Norm{y}_{\ell_\infty^N(\ell_2^k)} \,.
\]
The term $\Norm{\D' y - \D x}_{\ell_\infty^N(\ell_\infty^N)} \leq \Norm{\D' y - \D' Mx}_{\ell_\infty^N(\ell_\infty^N)} + \Norm{\D' Mx - \D x}_{\ell_\infty^N(\ell_\infty^N)}$
is estimated by
\[
    \Norm{\D' y - \D Mx}_{\ell_\infty^N(\ell_\infty^N)} &=& \max_{i,j} \biggl| \|y_i-y_j\|_{\ell_2^k} - \|Mx_i - Mx_j\|_{\ell_2^k} \biggr| \\
              &\leq& \max_{i,j} \|y_i-Mx_i\|_{\ell_2^k} + \|y_j - Mx_j\|_{\ell_2^k}
              \leq 2\E(t) \,,
\]
and, using the assumption~\eqref{Cont-JL},
\[
   \Norm{\D' Mx - \D x}_{\ell_\infty^N(\ell_\infty^N)} = \max_{i,j} \biggl| \|Mx_i-Mx_j\|_{\ell_2^k} - \|x_i - x_j\|_{\ell_2^d} \biggr|
                \leq \eps \max_{i,j} \|x_i - x_j\|_{\ell_2^k} = \eps \Norm{\D x}_{\ell_\infty^N(\ell_\infty^N)} \,.
\]
Finally, by the a priori estimate~\eqref{alpha} for $\Norm{\D x}_{\ell_\infty^N(\ell_\infty^N)}$
and~\eqref{beta} for $\Norm{y}_{\ell_\infty^N(\ell_2^d)}$, we obtain
\[
    \tot{}{t} e_i &\leq&  L \Norm{M} (2\E(t) + \eps\alpha) + L' \E(t) + L'' \beta (2\E(t) + \eps\alpha) \\
          &=& (2L\Norm{M} + 2\beta L'' + L')\E(t) + \eps\alpha(L\Norm{M} + L''\beta) \,.
\]

Now, let us split the interval $[0,T)$ into a union of finite disjoint intervals
$I_j=[t_{j-1},t_j)$, $j=1,\dots,K$
for a suitable $K\in\mathbb{N}$, such that $\E(t) = e_{i(j)}(t)$ for $t\in I_j$.
Consequently, on every $I_j$ we have
\[
    \tot{}{t} \E(t) = \tot{}{t} e_{i(j)}(t) \leq (2L\Norm{M} + 2\beta L'' + L')\E(t) + \eps\alpha(L\Norm{M} + L''\beta)  \,,
\]
and the Gronwall lemma yields
\[
    \E(t) \leq \left[\eps\alpha(L\Norm{M} + L''\beta)(t-t_{j-1}) +\E(t_{j-1})\right]
          \exp\left((2L\Norm{M} + 2\beta L'' + L')(t-t_{j-1})\right)
\]
for $t\in [t_{j-1},t_j)$.
A concatenation of these estimates over the intervals $I_j$ leads finally to the expected error estimate
\[
    \E(t) \leq \eps\alpha t (L\Norm{M} + L''\beta)
          \exp\left[(2L\Norm{M} + 2\beta L'' + L') t\right]  \,.
\]
\end{proof}

\subsection{A continuous Johnson-Lindenstrauss Lemma}
Let us now go through the assumptions we made in the formulation of Theorem~\ref{Thm:Cont}
and discuss how they restrict the validity and applicability of the result.
First of all, let us mention that~\eqref{alpha} and~\eqref{beta}
can be easily proven to hold for locally Lipschitz right-hand sides $f_i$ and $f_{ij}$ on finite time intervals.
Obviously, the critical point for the applicability of Theorem \ref{Thm:Cont} is the question
how to find a matrix $M$ satisfying the condition~\eqref{Cont-JL}, i.e.,
being a quasi-isometry along the trajectory solution $x(t)$ for {\it every} $t\in[0,T]$.
The answer is provided by the following generalization of the Johnson-Lindenstrauss Lemma
(Lemma~\ref{JLlem}) for rectifiable $\mathcal C^1$-curves, by a suitable continuity argument. 
Let us stress that our approach resembles the ``sampling and $\epsilon$-net'' argument in \cite{BDDW,BW09,W08} for the extension of
the quasi-isometry property of Johnson-Lindenstrauss embeddings to smooth Riemmanian manifolds. From this point of view the following result can be viewed
as a specification of the work  \cite{BW09,W08}.
\\
We first prove an auxiliary technical result:

\begin{lemma}\label{lem:neigh}
Let $0<\varepsilon<\varepsilon'<1$, $a\in\R^d$ and let $M:\R^d\to \R^k$ be a linear mapping such that
$$
(1-\varepsilon)\|a\|_{\ell_2^d}\le \|M a\|_{\ell_2^k}\le
(1+\varepsilon)\|a\|_{\ell_2^d}.
$$
Let $x\in\R^d$ satisfy 
\begin{equation}\label{eq:tau}
\|a-x\|\le \frac{(\varepsilon'-\varepsilon)\|a\|_{\ell_2^d}}{\|M\|+1+\varepsilon'}.
\end{equation}
Then
\begin{equation}\label{eq:neigh}
(1-\varepsilon')\|x\|_{\ell_2^d}\le \|Mx \|_{\ell_2^k}\le
(1+\varepsilon')\|x\|_{\ell_2^d}.
\end{equation}
\end{lemma}
\begin{proof} If $a=0$, the statement is trivial. If $a\not=0$, we denote the right-hand side of \eqref{eq:tau} by $\tau>0$ and
estimate by the triangle inequality
\begin{align*}
\frac{\|Mx\|_{\ell_2^k}}{\|x\|_{\ell_2^d}}&=\frac{\|M(x-a)+Ma \|_{\ell_2^k}}{\|x-a+a\|_{\ell_2^d}}\le
\frac{\|M\|\cdot\|x-a\|_{\ell_2^d}+(1+\varepsilon)\|a\|_{\ell_2^d}}{\|a\|_{\ell_2^d}-\|x-a\|_{\ell_2^d}}\\
&\le \frac{\|M\|\cdot\tau+(1+\varepsilon)\|a\|_{\ell_2^d}}{\|a\|_{\ell_2^d}-\tau}\le 1+\varepsilon' \,.
\end{align*}
A similar chain of inequalities holds for the estimate from below.
\end{proof}

Now we are ready to establish a continuous version of Lemma~\ref{JLlem}.

\begin{theorem}\label{thmcontJL}
Let $\varphi:[0,1]\to \R^d$ be a ${\mathcal C}^1$ curve. Let $0<\varepsilon<\varepsilon'<1$,
$$
\gamma:=\max_{\xi\in[0,1]}\frac{\|\varphi'(\xi)\|_{\ell_2^d}}{\|\varphi(\xi)\|_{\ell_2^d}}<\infty
\quad \text{and}\quad
\mathcal N\ge (\sqrt d+2)\cdot \frac{\gamma}{\varepsilon'-\varepsilon}.
$$
Let $k$ be such that a randomly chosen (and properly normalized) projector $M$ satisfies
the statement of the Johnson-Lindenstrauss Lemma~\ref{JLlem} with $\varepsilon, d, k$ and $\mathcal N$ arbitrary points with high probability. Without loss of generality we assume that $\|M\|\le \sqrt{d/k}$ within the same probability (this is in fact the case, e.g., for the examples of Gaussian and Bernoulli random matrices reported in Section \ref{Sec:Euler}).

Then 
\begin{equation}\label{eq:JLphi2}
(1-\varepsilon')\|\varphi(t)\|_{\ell_2^d}\le \|M\varphi(t)\|_{\ell_2^k}\le
(1+\varepsilon')\|\varphi(t)\|_{\ell_2^d},\mbox{ for all } t\in[0,1]
\end{equation}
holds with the same probability.
\end{theorem}

\begin{proof}
Let $t_i=i/\mathcal N$, $i=0,\dots,\mathcal N$  and put
$$
T_i:={\rm arg \ max}_{\xi\in [t_i,t_{i+1}]} \|\varphi'(\xi)\|_{\ell_2^d},\quad i=0,\dots,\mathcal N-1.
$$
Let $M:\R^d\to \R^k$ be the randomly chosen and normalized projector (see Lemma \ref{JLlem}). Hence $\|M\|\le \sqrt{d/k}$
and
\begin{equation}\label{eq:JL1}
(1-\varepsilon')\|\varphi(T_i)\|_{\ell_2^d}\le \|M(\varphi(T_i))\|_{\ell_2^k}\le
(1+\varepsilon')\|\varphi(T_i)\|_{\ell_2^d},\qquad i=1,\dots,\mathcal N
\end{equation}
with high probability. We show that \eqref{eq:JLphi2} holds with (at least) the same probability.

This follows easily from \eqref{eq:JL1} and the following estimate, which holds for every $t\in[t_i,t_{i+1}]$,
\begin{align*}
\|\varphi(t)-\varphi(T_i)\|_{\ell_2^d}&\le \int_{t}^{T_i} \|\varphi'(s)\|_{\ell_2^d}ds\le \frac{\|\varphi'(T_i)\|_{\ell_2^d}}{\mathcal N}\le
\frac{\|\varphi'(T_i)\|_{\ell_2^d}(\varepsilon'-\varepsilon)}{\gamma(\sqrt d+2)}\\
&\le \frac{\|\varphi(T_i)\|_{\ell_2^d}(\varepsilon'-\varepsilon)}{\sqrt d+2}
\le \frac{\|\varphi(T_i)\|_{\ell_2^d}(\varepsilon'-\varepsilon)}{\|M\|+1+\varepsilon'}.
\end{align*}
The proof is then finished by a straightforward application of Lemma \ref{lem:neigh}.
\end{proof}
\begin{remark}
We show now that the condition $$
\gamma:=\max_{\xi\in[0,1]}\frac{\|\varphi'(\xi)\|_{\ell_2^d}}{\|\varphi(\xi)\|_{\ell_2^d}}<\infty$$
is necessary, hence it is a restriction to the type of curves one can quasi-isometrically project. Let $d\ge 3$.
It is known that there is a continuous curve $\varphi:[0,1]\to [0,1]^{d-1}$,
such that $\varphi([0,1])=[0,1]^{d-1}$, i.e., $\varphi$ goes \emph{onto} $[0,1]^{d-1}.$ 
The construction of such a \emph{space-filling} curve goes back to Peano and Hilbert. After a composition with
suitable dilations and $d$-dimensional spherical coordinates we observe that there is also
a \emph{surjective} continuous curve $\phi:[0,1]\to {\mathbb S}^{d-1}$, where ${\mathbb S}^{d-1}$
denotes the $\ell_2^d$ unit sphere in $\R^d.$

As $M$ was supposed to be a projection, \eqref{eq:JLphi2} cannot hold for all $t$'s with
$\varphi(t)\in {\rm ker}\ M\not = \emptyset.$
\end{remark}

Obviously, the key condition for applicability of Theorem \ref{thmcontJL} for finding
a projection matrix $M$ satisfying~\eqref{Cont-JL} is that
\(  \label{Cond-contJL}
   \sup_{t\in[0,T]} \max_{i,j} \frac{\|\dot x_i - \dot x_j\|_{\ell_2^d}}{\|x_i-x_j\|_{\ell_2^d}} \leq \gamma <\infty \,.
\)
This condition is, for instance, trivially satisfied when the right-hand sides $f_i$'s and $f_{ij}$'s
have the following Lipschitz continuity:
\[   
   \|f_i(\D x)  - f_j(\D x)\|_{\ell_2^d} \leq L''' \|x_i-x_j\|_{\ell_2^d} \qquad\mbox{for all } i,j=1,\dots,N \,,\\
   |f_{i,k}(\D x)  - f_{j,k}(\D x)| \leq L'''' \|x_i-x_j\|_{\ell_2^d} \qquad\mbox{for all } i,j,k=1,\dots,N.
\]
We will show in the examples below how condition~\eqref{Cond-contJL} is verified
in cases of dynamical systems modeling standard social mechanisms of
{\it attraction, repulsion, aggregation} and {\it alignment}.

\subsection{Applicability to relevant examples of dynamical systems describing social dynamics}

\begin{figure}[ht]
\[\begin{array}{ccc}
   \centering
     \includegraphics[width=4cm]{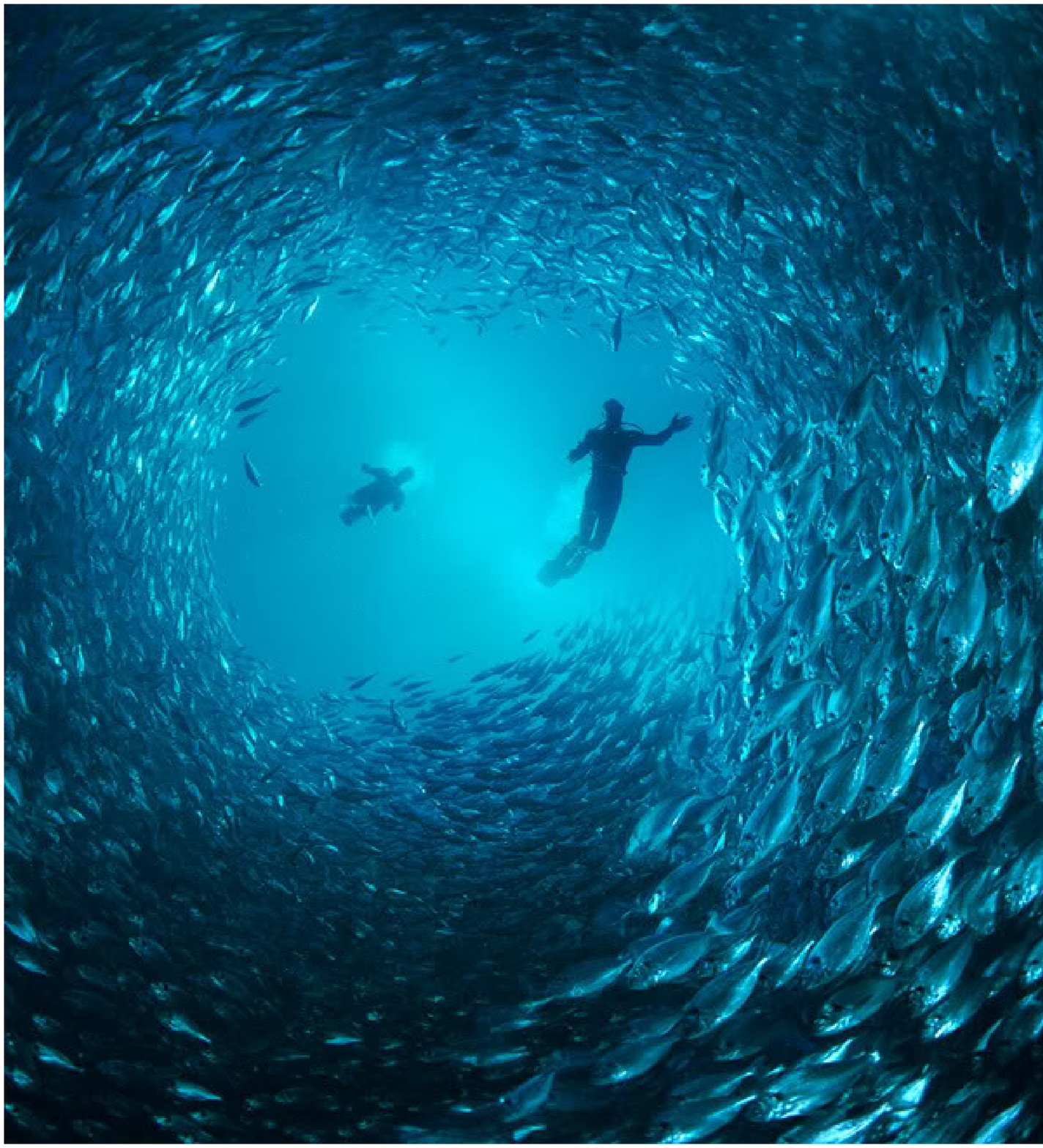} & &
 \includegraphics[width=8cm]{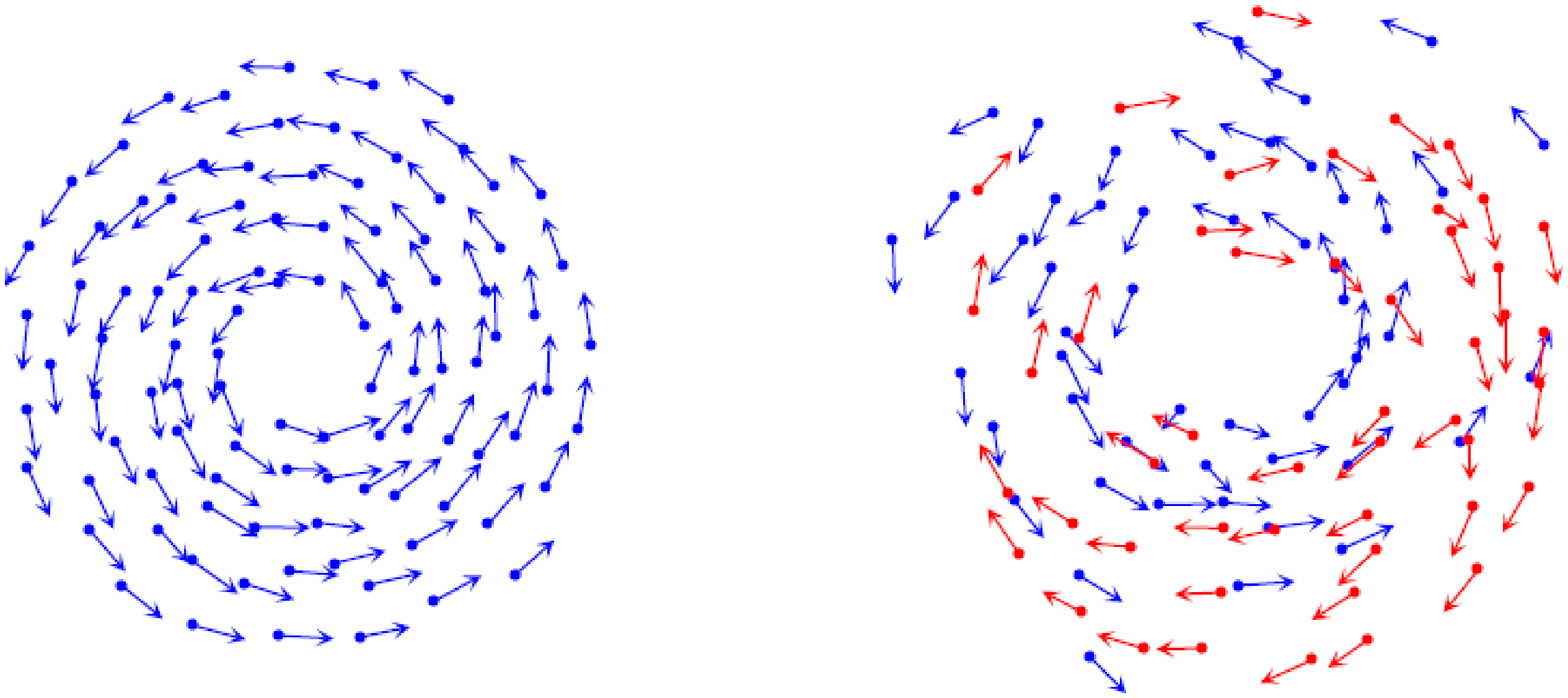}
\end{array}\]
\caption{Mills in nature and in models}
\label{fig:mills}
\end{figure}

In this section we show the applicability of our dimensionality reduction theory to well-known dynamical systems
driven by ``social forces'' of {\it alignment, attraction, repulsion, aggregation, and self-drive}.
Although these models were proposed as descriptions of {\it group motion in physical space},
the fundamental social effects can be considered as building blocks in the more abstract context
of many-parameter social dynamics. 
It has been shown \cite{CFTV,DOrsogna} that these models are able to produce meaningful {\it patterns},
for instance {\it mills} in two spatial dimensions (see Figure~\ref{fig:mills}),
reproducing the behavior of certain biological species.
However, we should expect that in higher dimension the possible patterns produced by the combination of fundamental effects can be much more complex.

\subsubsection{The Cucker-Smale system (alignment effect)}\label{sec:cs}
As shown in Section~\ref{Sec:Euler}, the Cucker and Smale flocking model~\eqref{CS1}--\eqref{CS2}
is of the type~\eqref{eq:dyn}, satisfies the Lipschitz continuity assumptions~\eqref{eq:condf1}--\eqref{eq:condf3},
and it is bounded at finite time, as already discussed in Section \ref{Subs:UniformCS}.
Therefore, to meet all the assumptions of Theorem~\ref{Thm:Cont},
we only need to check that it also satisfies the condition~\eqref{Cond-contJL}.
However, for this we need to consider a slightly different framework than in Section~\ref{Subs:UniformCS};
instead of considering the $2N$ $d$-dimensional variables ($N$ position variables
and $N$ velocity variables), we need to arrange the model as $N$ variables in $\R^{2d}$,
each variable consisting of the position part (first $d$ entries) and of the velocity part
(the other $d$ entries).
We have then
\[
    \|\dot x_i-\dot x_j\|_{\ell_2^d} + \|\dot v_i -\dot v_j\|_{\ell_2^d} &\leq& \|v_i-v_j\|_{\ell_2^d}
          + \frac1N \sum_{k=1}^N \bigl |g(\|x_i-x_k\|_{\ell_2^d}) - g(\|x_j-x_k\|_{\ell_2^d})\bigr | \|v_k\|_{\ell_2^d}   \\
      &\leq& \|v_i-v_j\|_{\ell_2^d} + \frac{\Norm{g}_{Lip}}{N} \sum_{k=1}^N \bigl| \|x_i-x_k\|_{\ell_2^d} - \|x_j-x_k\|_{\ell_2^d}\bigr| \|v_k\|_{\ell_2^d} \\
      &\leq& \|v_i-v_j\|_{\ell_2^d} + \frac{\Norm{g}_{Lip}}{N} \left( \sum_{k=1}^N \|v_k\|_{\ell_2^d}\right) \|x_i-x_j\|_{\ell_2^d} \\
      &\leq& \|v_i-v_j\|_{\ell_2^d} + c \|x_i-x_j\|_{\ell_2^d} \,,
\]
for a suitable constant $c$ depending on the initial data.
We used here the a-priori boundedness
of the term $\frac{1}{N} \left( \sum_{k=1}^N \|v_k\|_{\ell_2^d}\right)$,
see~\cite{CS2} or~\cite{TH} for details.
Consequently, we can satisfy~\eqref{Cond-contJL} with $\gamma=\max(1,c)$.

\subsubsection{Second order dynamic model with self-propulsion and pairwise interactions (self-drive, attraction, and repulsion effects)}\label{sec:dors}
Another practically relevant model which fits into the class given by~\eqref{eq:dyn}
is a second order dynamic model with self-propulsion and pairwise interactions,~\cite{LRC00,DOrsogna}:
\(
    \dot x_i &=& v_i \,, \label{DOrsogna1}\\
    \dot v_i &=& (a-b\|v_i\|_{\ell_2^d}^2)v_i - \frac{1}{N} \sum_{j\neq i} \grad_{x_i} U(\|x_i-x_j\|_{\ell_2^d}) \,,\qquad i=1,\dots,N,
      \label{DOrsogna2}
\)
where $a$ and $b$ are positive constants and $U:[0,\infty)\to\R$ is a smooth potential.
We denote $u(s) = U'(s)/s$ and assume that $u$ is a bounded, Lipschitz continuous function.
We again arrange the model as a system of $N$ variables in $\R^{2d}$,
each variable consisting of the position part (first $d$ entries) and of the velocity part
(the other $d$ entries). Consequently, the model can be put into a form compliant with~\eqref{eq:dyn} as follows:
\[
    \dot x_i &=& \sum_{j=1}^N f^{xv}_{ij} v_j \,,\\
    \dot v_i &=& \sum_{j=1}^N f^{vv}_{ij}(\D v) v_j + \sum_{j=1}^N f^{vx}_{ij}(\D x) x_j\,,\\
\]
with $f^{xv}_{ij} = \delta_{ij}$,
$f^{vx}_{ii}(\D x) = - \frac{1}{N} \sum_{j\neq i} u(\|x_i-x_j\|_{\ell_2^d})$ and
$f^{vx}_{ij}(\D x) = \frac{1}{N}  u(\|x_i-x_j\|_{\ell_2^d})$ for $i\neq j$.
Moreover, we may set  $f^{vv}_{ij}(\D v) = \delta_{ij} (a-b\|v_i\|_{\ell_2^d}^2)$
by introducing an auxiliary, noninfluential constant zero particle $(x_0,v_0)=(0,0)$
with null dynamics, i.e., $f_0^{*\star}=0$ and $f_{0j}^{*\star}=0$, where $*,\star \in \{x,v\}$.
Then,~\eqref{eq:condf1} is void, while~\eqref{eq:condf2} is satisfied by
\begin{align*}
    \max_i \sum_j (|f^{xv}_{ij}(\D x,\D v)|&+
    |f^{vx}_{ij}(\D x,\D v)|+|f^{vv}_{ij}(\D x,\D v)|)\\
    & \leq 1 + a + b\max_i \|v_i\|_{\ell_2^d}^2 + 2 \Norm{u}_{L_\infty} \leq L' \,,
\end{align*}
since the theory provides an apriori bound on $\beta_v:=\sup_{t\in[0,T]} \max_i \|v_i\|_{\ell_2^d}$, see~\cite{DOrsogna}.
Condition~\eqref{eq:condf3} for $f^{xv}_{ij}$ is void, while for $f^{vv}_{ij}$ it is satisfied by
\begin{align*}
    \max_i \sum_j \left| f_{ij}^{vv}(\D v) - f_{ij}^{vv}(\D w) \right|
      &\leq b\max_i \left| \|v_i\|_{\ell_2^d}^2 - \|w_i\|_{\ell_2^d}^2 \right| \\
&\leq b \max_i\left(\|v_i\|_{\ell_2^d}+\|w_i\|_{\ell_2^d}\right) \|v_i-w_i\|_{\ell_2^d}\\
      &\leq L''\Norm{\D v - \D w}_{\ell_\infty^N(\ell_\infty^N)} \,,
\end{align*}
where we again use the apriori boundedness of $\beta_v$.
For $f^{vx}_{ij}$ is~\eqref{eq:condf3} satisfied by
\[
    \max_i \sum_j \left| f_{ij}^{vx}(\D x) - f_{ij}^{vx}(\D y) \right|
      &\leq& \max_i \frac{2}{N} \sum_{j\neq i} \left| u(\|x_i-x_j\|_{\ell_2^d}) - u(\|y_i-y_j\|_{\ell_2^d})\right| \\
      &\leq& \max_i \frac{2}{N} \Norm{u}_{{\rm Lip}} \sum_{j\neq i} \left| \|x_i-x_j\|_{\ell_2^d} - \|y_i-y_j\|_{\ell_2^d} \right| \\
      &\leq& 2 \Norm{u}_{{\rm Lip}} \Norm{\D x - \D y}_{\ell_\infty^N(\ell_\infty^N)}.
\]
Finally, it can be easily checked that condition~\eqref{Cond-contJL} is satisfied by
\[
    \|\dot x_i - \dot x_j\|_{\ell_2^d} + \|\dot v_i - \dot v_j\|_{\ell_2^d} \leq
         (1+a+3b\beta_v^2) \|v_i-v_j\|_{\ell_2^d} + \left( \Norm{u}_{L_\infty} + 2 \beta_x \Norm{u}_{{\rm Lip}} \right) \|x_i-x_j\|_{\ell_2^d} \,,
\]
where $\beta_x:=\sup_{t\in[0,T]} \max_i \|x_i\|_{\ell_2^d}$.
We notice that also this model is bounded at finite time as shown in \cite[Theorem 3.10 (formula (22))]{CCR}, and therefore for any
fixed horizon time $T$, there is a constant $\alpha=\alpha(T)>0$ such that \eqref{alpha1} and \eqref{alpha} hold.
In the paper \cite{DOrsogna}
it is shown that this model tends to produce patterns of different quality, in particular
mills, double mills, and translating crystalline flocks (see also Figure \ref{fig:mills}).
These patterns were further studied in~\cite{cop}.
Starting from the Liouville equation for the many-body problem the authors derive
the corresponding kinetic and macroscopic hydrodynamic equations.
The kinetic theory approach leads to the identification of macroscopic structures
otherwise not recognized as solutions of the hydrodynamic equations,
such as double mills of two superimposed flows. The authors found conditions
allowing for the existence of such solutions and compared them to the case of single mills.
In~\cite{combc} the authors utilize the methods of classical statistical mechanics
to connect the individual-based models of the type~\eqref{DOrsogna1}--\eqref{DOrsogna2} to
their continuum formulations and determine criteria for the validity of the latter.
They show that H-stability of the interaction potential plays a fundamental role in determining both
the validity of the continuum approximation and the nature of the aggregation state
transitions. They perform a linear stability analysis of the continuum model and compare
the results to the simulations of the individual-based one.

Without entering into further details, let us stress that mills and double
mills are uniformly bounded in time (and stable). Hence in these cases,
we can assume that actually the constant $\alpha$ is again bounded for all times. Moreover,
when the dynamics converges to a translating crystalline flocks, we may reason in a similar way as done 
for the Cucker-Smale model (although in this case the pattern in unstable).


\subsection{Recovery of the dynamics in high dimension from multiple simulations in low dimension}

The main message of Theorem \ref{Thm:Cont} is that, under suitable assumptions on the governing functions $f_i, f_{ij}$, the trajectory of the solution $y(t)$ of the {\it projected} dynamical system \eqref{lowdimsys} is at an $\varepsilon$ error from the trajectory of the {\it projection} of the solution $x(t)$ of the dynamical system  \eqref{Cont1}-\eqref{Cont2}, i.e.,
\begin{equation}\label{approxk}
y_i (t) \approx M x_i(t) \mbox { or, more precisely, } \| Mx_i(t) - y_i(t) \|_{\ell_2^k} \leq C(t) \varepsilon, \quad t \in [0,T].
\end{equation}

We wonder whether this approximation property can allow us to ``learn'' properties of the original trajectory $x(t)$ in high dimension. 

 \subsubsection{Optimal information recovery of high-dimensional trajectory from low-dimensional projections}\label{CSsec}

In this section we would like to address the following two fundamental questions:
\begin{itemize}
\item[(i)] Can we quantify the best possible information of the high-dimensional trajectory one can recover from one or more projections in
lower dimension?
\item[(ii)] Is there any practical method which performs an optimal recovery?  
\end{itemize}
The first question was implicitly addressed already in the 70's  by Kashin and later by Garnaev and Gluskin \cite{ka77,gagl84}, as one can put in relationship the optimal recovery from linear measurements 
with Gelfand width of $\ell_p$-balls, see for instance \cite{codade09}. It was only with the development of the theory of {\it compressed sensing} \cite{carota06,do06-2} that an answer to the second question was provided, showing that $\ell_1$-minimization actually performs an {\it optimal recovery} of vectors in high dimension from random linear projections to low dimension. We address the reader to \cite[Section 3.9]{FR}  for further details.
In the following we concisely recall the theory of compressed sensing and we apply it to estimate the optimal information error
in recovering the trajectories in high dimension from lower dimensional simulations. 

Again a central role here is played by (random) matrices with the so-called {\it Restricted Isometry Property} RIP, cf. \cite{ca08}.

\begin{definition}[Restricted Isometry Property]
A $k\times d$ matrix $M$ is said to have the Restricted Isometry Property of order $K\leq d$ and
level $\delta\in (0,1)$ if
$$
(1-\delta)\|x\|^2_{\ell_2^d} \leq \|M x\|^2_{\ell_2^k} \leq (1+\delta) \|x\|^2_{\ell_2^d}
$$
for all $K$-sparse $x \in \Sigma_K=\{z \in \mathbb R^d: \#\supp(z) \leq K\}$.
\end{definition}

Both the typical matrices used in Johnson-Lindenstrauss embeddings (cf. Lemma \ref{JLlem}) and matrices
with RIP used in compressed sensing are usually generated at random. It was observed by \cite{BDDW} and \cite{KW},
that there is an intimate connection between these two notions. A simple reformulation
of the arguments of \cite{BDDW} yields the following.

\begin{theorem}[Baraniuk, Davenport, DeVore, and Wakin]\label{BDDWthm}
Let $M$ be a $k\times d$ matrix drawn at random which satisfies
$$
(1-\delta/2)\|x\|^2_{\ell_2^d} \leq \|M x\|^2_{\ell_2^k} \leq (1+\delta/2) \|x\|^2_{\ell_2^d},\quad x\in {\mathcal P}
$$
for every set ${\mathcal P}\subset \R^d$ with $\#{\mathcal P}\le \bigl(\frac{12ed}{\delta K}\bigr)^K$ with probability $0<\nu<1$.
Then $M$ satisfies the Restricted Isometry Property of order $K$ and level $\delta/3$ with probability at least equal to $\nu$. 
\end{theorem}

Combined with several rather elementary constructions of Johnson-Lindenstrauss embedding matrices available in literature, cf. \cite{A}
and \cite{DG}, this result provides a simple construction of RIP matrices. The converse direction, namely the way from RIP matrices
to matrices suitable for Johnson-Lindenstrauss embedding was discovered only recently in \cite{KW}.

\begin{theorem}[Krahmer and Ward]\label{KWthm}
Fix $\eta > 0$ and $\varepsilon> 0$, and consider a finite set $\mathcal P\subset \mathbb R^d$ of cardinality $|\mathcal P| = \mathcal N$. Set
$K \geq 40 \log \frac{4 \mathcal N}{\eta}$, and suppose that the $k \times d$ matrix $\tilde M$  satisfies  the Restricted Isometry Property of order $K$ and
level $\delta \leq \varepsilon/4$. Let $\xi \in \mathbb R^d$ be a Rademacher sequence, i.e., uniformly distributed on $\{-1, 1\}^d$ . Then with
probability exceeding $1 -\eta$,
\begin{equation*}
(1-\varepsilon)\|x\|_{\ell_2^d}^2\le \|M x\|_{\ell_2^k}^2\le (1+\varepsilon)\|x\|_{\ell_2^d}^2.
\end{equation*}
uniformly for all $x \in \mathcal P$, where $M:= \tilde M \operatorname{diag}(\xi)$, where $\operatorname{diag}(\xi)$
is a $d\times d$ diagonal matrix with $\xi$ on the diagonal.
\end{theorem}

We refer to \cite{ra10} for additional  details.

\begin{remark}
Notice that $M$ as constructed in Theorem \ref{KWthm} is both a Johnson-Lindenstrauss embedding and a matrix with RIP, because
\begin{align*}
(1-\delta)\|x\|_{\ell_2^d}^2 &= (1-\delta)\|\operatorname{diag}(\xi) x\|_{\ell_2^d}^2 \leq  \|\underbrace{\tilde M \operatorname{diag}(\xi)}_{:=M} x\|_{\ell_2^k}^2 \\
&\leq (1+\delta) \|\operatorname{diag}(\xi) x\|_{\ell_2^d}^2 = (1+\delta) \| x\|_{\ell_2^d}^2.
\end{align*}
The matrices considered in Section \ref{Sec:Euler} satisfy with high probability the RIP with
$$
K = \mathcal O \left ( \frac{k}{1+ \log( d/k)} \right).
$$

\end{remark}

Equipped with the notion of RIP matrices we may state the main result of the theory of compressed sensing, as appearing in \cite{fo09}, which we shall
use for the recovery of the dynamical system in $\R^d.$
\begin{theorem}\label{thmStrongRIP}
Assume that the matrix $M \in \mathbb R^{k \times d}$ has the RIP of order $2 K$ and level 
$$
\delta_{2 K} < \frac{2}{3+\sqrt{7/4}} \approx 0.4627. 
$$
Then the following holds for all $x \in \mathbb R^d$. Let the low-dimensional approximation $y = M x + \eta$ be given
with $\|\eta\|_{\ell_2^k}\leq C \varepsilon$. Let $x^\#$ be the solution of
\begin{equation}\label{P1:eta}
\min_{z \in \mathbb R^d} \|z\|_{\ell_1^d} \quad \mbox{ subject to } \|M z - y\|_{\ell_2^k} \leq \|\eta\|_{\ell_2^k}.
\end{equation}
Then 
\[
\|x - x^\#\|_{\ell_2^d} \leq C_1 \varepsilon + C_2 \frac{\sigma_K(x)_{\ell_1^d}}{\sqrt{K}}
\]
for some constants $C_1,C_2 > 0$ that depend only on $\delta_{2K}$, and $\sigma_K(x)_{\ell_1^d} = \inf_{z: \#\supp(z) \leq K} \|z -x\|_{\ell_1^d}$ is the best-$K$-term approximation error in $\ell_1^d$.
\end{theorem}
\\ 
This result says that provided the stability relationship \eqref{approxk}, we can approximate the individual trajectories $x_i(t)$, for each $t \in [0,T]$ fixed, 
by a vector $x^\#_i(t)$ solution of an optimization problem of the type \eqref{P1:eta}, and the accuracy of the approximation depends on the best-$K$-term approximation 
error $\sigma_K(x_i(t))_{\ell_1^d}$. Actually, the results in \cite{carota06,do06-2} in connection with \cite{codade09,ka77,gagl84}, state also that this is asymptotically the best one can hope for.
One possibility to improve the recovery error is to increase the dimension $k$ (leading to a smaller distortion parameter $\varepsilon>0$ in the Johnson-Lindenstrauss
embedding). But we would like to explore another possibility, namely
projecting and simulating {\it in parallel and independently} the dynamical system $L$-times in the lower dimension $k$
\begin{equation}\label{eq:paral1}
\dot y^\ell_i = M^\ell f_i(\D' y^\ell) + \sum_{j=1}^N f_{ij}(\D' y^\ell) y_j^\ell\,,\qquad y_i^\ell(0) = M^\ell x_i^0 \,,\quad \ell=1,\dots, L.
\end{equation}
Let us give a brief overview of the corresponding error estimates. The number of points needed in each of the cases is
${\mathcal N}\approx N\times n_0$, where $N$ is the number of agents and $n_0=T/h$ is the number of iterations.
\begin{itemize}
\item We perform 1 projection and simulation in $\R^k$: Then $\varepsilon={\mathcal O}\Bigl(\sqrt{\frac{\log {\mathcal N}}{k}}\Bigr)$, $K=\mathcal O \left ( \frac{k}{1+ \log( d/k)} \right)$
and an application of Theorem \ref{thmStrongRIP} leads to
\begin{equation}\label{eq:paral2}
\|x_i(t) - x_i^\#(t)\|_{\ell_2^d} \le C'(t)\left( \sqrt{\frac{\log {\mathcal N}}{k}} + \frac{\sigma_K(x_i(t))_{\ell_1^d}}{\sqrt{K}}\right).
\end{equation}
Here, $C'(t)$ combines both the constants from Theorem \ref{thmStrongRIP} and the time-dependent $C(t)$ from \eqref{approxk}.
So, to reach the precision of order $C'(t)\epsilon>0$, we have to choose $k\in\N$ large enough, such that
$\sqrt{\frac{\log {\mathcal N}}{k}}\le \epsilon$ and $\frac{\sigma_K(x_i(t))_{\ell_1^d}}{\sqrt{K}}\le \epsilon$. We then need
$k\times N^2$ operations to evaluate the adjacency matrix.

\item We perform 1 projection and simulation in $\R^{L\times k}$: Then $\varepsilon'={\mathcal O}\Bigl(\sqrt{\frac{\log {\mathcal N}}{Lk}}\Bigr)$ and
$K'=\mathcal O \left ( \frac{Lk}{1+ \log(d/Lk)} \right)$ and an application of Theorem \ref{thmStrongRIP} leads to
\begin{equation}\label{eq:paral4}
\|x_i(t) - x_i^\#(t)\|_{\ell_2^d} \le C'(t)\left( \sqrt{\frac{\log {\mathcal N}}{Lk}} + \frac{\sigma_{K'}(x_i(t))_{\ell_1^d}}{\sqrt{K'}}\right).
\end{equation}

The given precision of order $C'(t)\epsilon>0$, may be then reached by choosing $k,L\in\N$ large enough, such that
$\sqrt{\frac{\log {\mathcal N}}{Lk}}\le \epsilon$ and
$\frac{\sigma_{K'}(x_i(t))_{\ell_1^d}}{\sqrt{K'}}\le \epsilon$.
We then need $Lk\times N^2$ operations to evaluate the adjacency matrix.

\item We perform $L$ independent and parallel projections and simulations in $\R^k$:
Then we assemble the following system corresponding to \eqref{eq:paral1}
$$
{\mathcal M}x=\left (
\begin{array}{l}
M^1 \\M^2\\\dots\\\dots \\ M^{L}
\end{array}
\right ) x_i = \left (
\begin{array}{l}
y^1_i \\y^2_i\\\dots\\\dots \\y^{L}_i
\end{array}
\right ) - \left (
\begin{array}{l}
\eta^1_i \\\eta^2_i\\\dots\\\dots \\\eta^{L}_i
\end{array}
\right ),
$$
where for all $\ell=1,\dots, L$ the matrices $M^\ell \in \mathbb R^{k \times d}$ are (let us say) random matrices with each entry generated
independently with respect to the properly normalized Gaussian distribution as described in Section \ref{Sec:Euler}.
Then ${\mathcal M}/\sqrt{L}$ is a $Lk\times d$ matrix with Restricted Isometry Property of order $K'=\mathcal O \left ( \frac{Lk}{1+ \log(d/Lk)} \right)$ and
level $\delta <  0.4627$. The initial distortion of each of the projections is still $\varepsilon={\mathcal O}\Bigl(\sqrt{\frac{\log {\mathcal N}}{k}}\Bigr)$.
Therefore, by applying  Theorem \ref{thmStrongRIP}, we can compute $x_i^\#(t)$ such that
\begin{equation}\label{eq:paral3}
\|x_i(t) - x^\#_i(t)\|_{\ell_2^d} \le C'(t)\left( \sqrt{\frac{\log {\mathcal N}}{k}} + \frac{\sigma_{K'}(x_i(t))_{\ell_1^d}}{\sqrt{K'}}\right).
\end{equation}
Notice that the computation of $x_i^\#(t)$ can also be performed in parallel, see, e.g., \cite{fo07}.
The larger is the number $L$ of projections we perform, the larger is $K'$ and the smaller is the second summand in \eqref{eq:paral3};
actually $\sigma_{K'}(x_i(t))_{\ell_1^d}$ vanishes for $K'\geq d$. Unfortunately, the parallelization can not help to reduce
the initial distortion $\varepsilon>0$. To reach again the precision of order $C'(t)\epsilon>0$, we have to choose $k\in\N$ large enough, 
such that
$\sqrt{\frac{\log {\mathcal N}}{k}}\le \epsilon$. Then we chose $L\ge 1$ large enough such that
$\frac{\sigma_{K'}(x_i(t))_{\ell_1^d}}{\sqrt{K'}}\le \epsilon$. We again need $k\times N^2$ operations to evaluate the adjacency matrix.
\end{itemize}

In all three cases, we obtain the estimate
\begin{equation}\label{eq:paral5}
\|x_i(t) - x_i^\#(t)\|_{\ell_2^d} \le C'(t)\left( \varepsilon + \frac{\sigma_{K}(x_i(t))_{\ell_1^d}}{\sqrt{K}}\right),
\end{equation}
where the corresponding values of $\varepsilon>0$ and $K$ together with the number of operations
needed to evaluate the adjacency matrix may be found in the following table.

\vskip.5cm
\begin{center}
\begin{tabular}{|c|c|c|c|}
\hline
&$\varepsilon$&$K$& number of operations\\
\hline
\vphantom{$\biggl($}1 projection into $\R^k$ & ${\mathcal O}\Bigl(\sqrt{\frac{\log {\mathcal N}}{k}}\Bigr)$ &
$\mathcal O \left ( \frac{k}{1+ \log( d/k)} \right)$ & $k\times N^2$\\
\hline
\vphantom{$\biggl($}1 projection into $\R^{L\times k}$ & ${\mathcal O}\Bigl(\sqrt{\frac{\log {\mathcal N}}{Lk}}\Bigr)$ &
$\mathcal O \left ( \frac{Lk}{1+ \log(d/Lk)} \right)$& $Lk\times N^2$\\
\hline
\vphantom{$\biggl($}$L$ projections into $\R^k$ & ${\mathcal O}\Bigl(\sqrt{\frac{\log {\mathcal N}}{k}}\Bigr)$ &
$\mathcal O \left ( \frac{Lk}{1+ \log(d/Lk)} \right)$ & $k\times N^2$\\
\hline
\end{tabular}
\end{center}

\subsubsection{Optimal recovery of trajectories on smooth manifolds}

In recent papers \cite{BW09,W08,IM11}, the concepts of compressed sensing and optimal recovery
were extended to vectors on smooth manifolds. These methods could become very useful
in our context if (for any reason) we would have an apriori knowledge that the trajectories $x_i(t)$
keep staying on or near such a smooth manifold. Actually this is the case, for instance in molecular dynamics,
where simulations, e.g. in the form of the coordinates of the atoms in a molecule as a function of time, lie on or near an intrinsically-low-dimensional set in the high-dimensional state space of the molecule, 
and geometric properties of such sets provide important information about the dynamics, or about how to build low-dimensional representations of such dynamics \cite{RZMC,ZRMC}.
In this case, by using appropriate recovery methods
as described in \cite{IM11}, we could recover high-dimensional vectors from very low dimensional
projections with much higher accuracy. However, this issue will be addressed in a following paper.

\subsection{Numerical experiments}\label{Sec:Numerics}
In this section we illustrate the practical use and performances of our projection method
for the Cucker-Smale system~\eqref{CS1}--\eqref{CS2}.

\subsubsection{Pattern formation detection in high dimension from lower dimensional projections}
As already mentioned, this system models the emergence of consensus in a group of
interacting agents, trying to align with their neighbors.
The qualitative behavior of its solutions is formulated by this well known result~\cite{CS1, CS2, TH}:

\begin{theorem}
Let $(x_i(t), v_i(t))$ be the solutions of~\eqref{CS1}--\eqref{CS2}.
Let us define the fluctuation of positions around the center of mass
$x_c(t) = \frac{1}{N} \sum_{i=1}^N x_i(t)$, and, resp.,
the fluctuation of the rate of change around its average $v_c(t) = \frac{1}{N} \sum_{i=1}^N v_i(t)$ as
\[
    \Lambda(t) = \frac1N \sum_{i=1}^N \|x_i(t) - x_c(t)\|_{\ell_2^d}^2 \,,\qquad
    \Gamma(t) = \frac1N \sum_{i=1}^N \|v_i(t) - v_c(t)\|_{\ell_2^d}^2 \,.
\]
Then if either $\beta\leq 1/2$ or the initial fluctuations $\Lambda(0)$ and $\Gamma(0)$
are small enough (see~\cite{CS1} for details),
then $\Gamma(t)\to 0$ as $t\to\infty$.
\end{theorem}

The phenomenon of $\Gamma(t)$ tending to zero as $t\to\infty$ is called \emph{flocking} or {\it emergence of consensus}.
If $\beta > 1/2$ and the initial fluctuations are not small,
it is not known whether a given initial configuration
will actually lead to flocking or not, and the only way to find out the possible formation of {\it consensus patterns}
is to perform  numerical simulations. However, these can be especially costly
if the number of agents $N$ and the dimension $d$ are large;
the algorithmic complexity of the calculation is $\O(d\times N^2)$.
Therefore, a significant reduction of the dimension $d$,
which can be achieved by our projection method,
would lead to a corresponding reduction of the computational cost.

\begin{figure}[ht]
\[\begin{array}{ccc}
   \centering
      \epsfig{figure={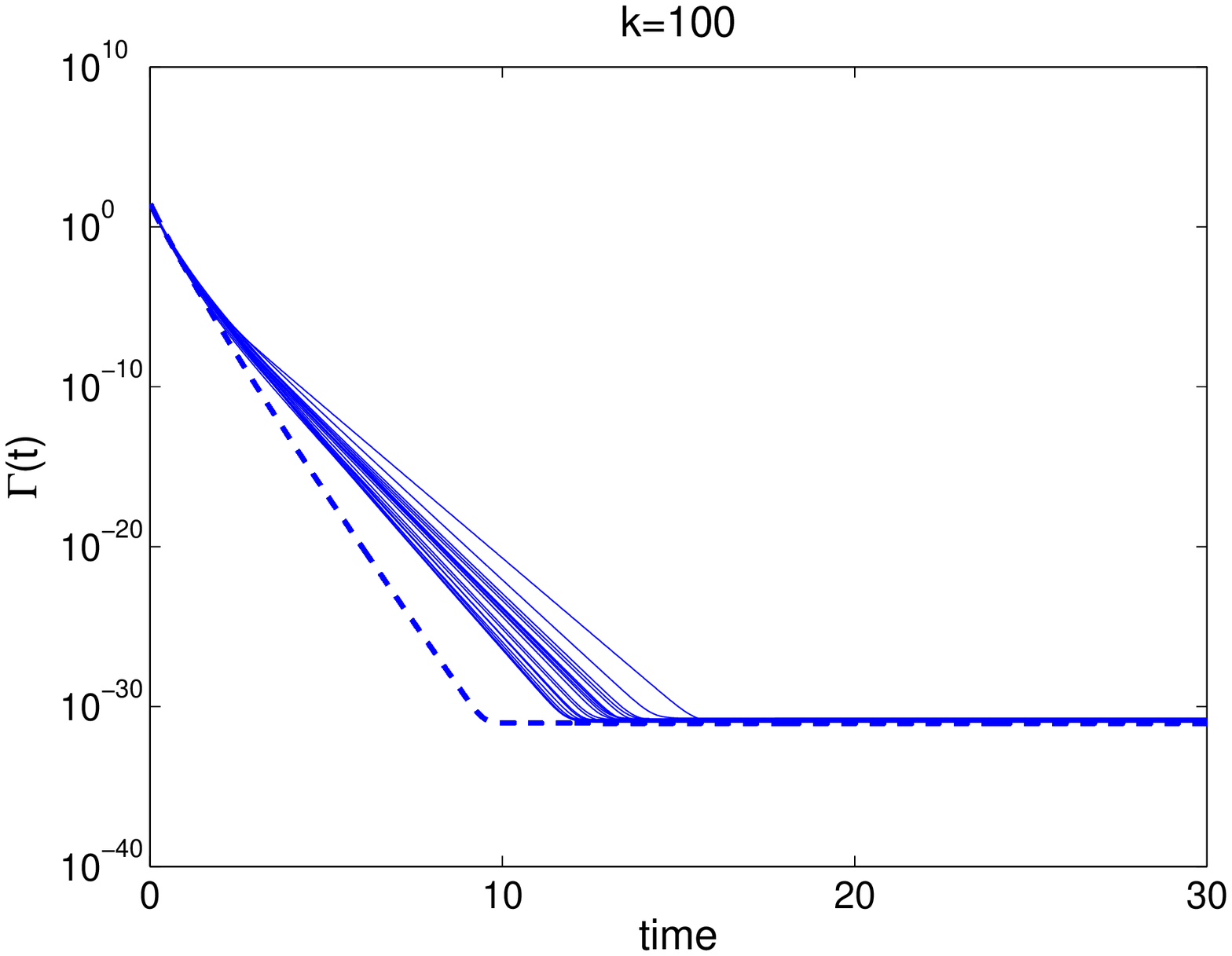}, width=0.5\linewidth} & &
      \epsfig{figure={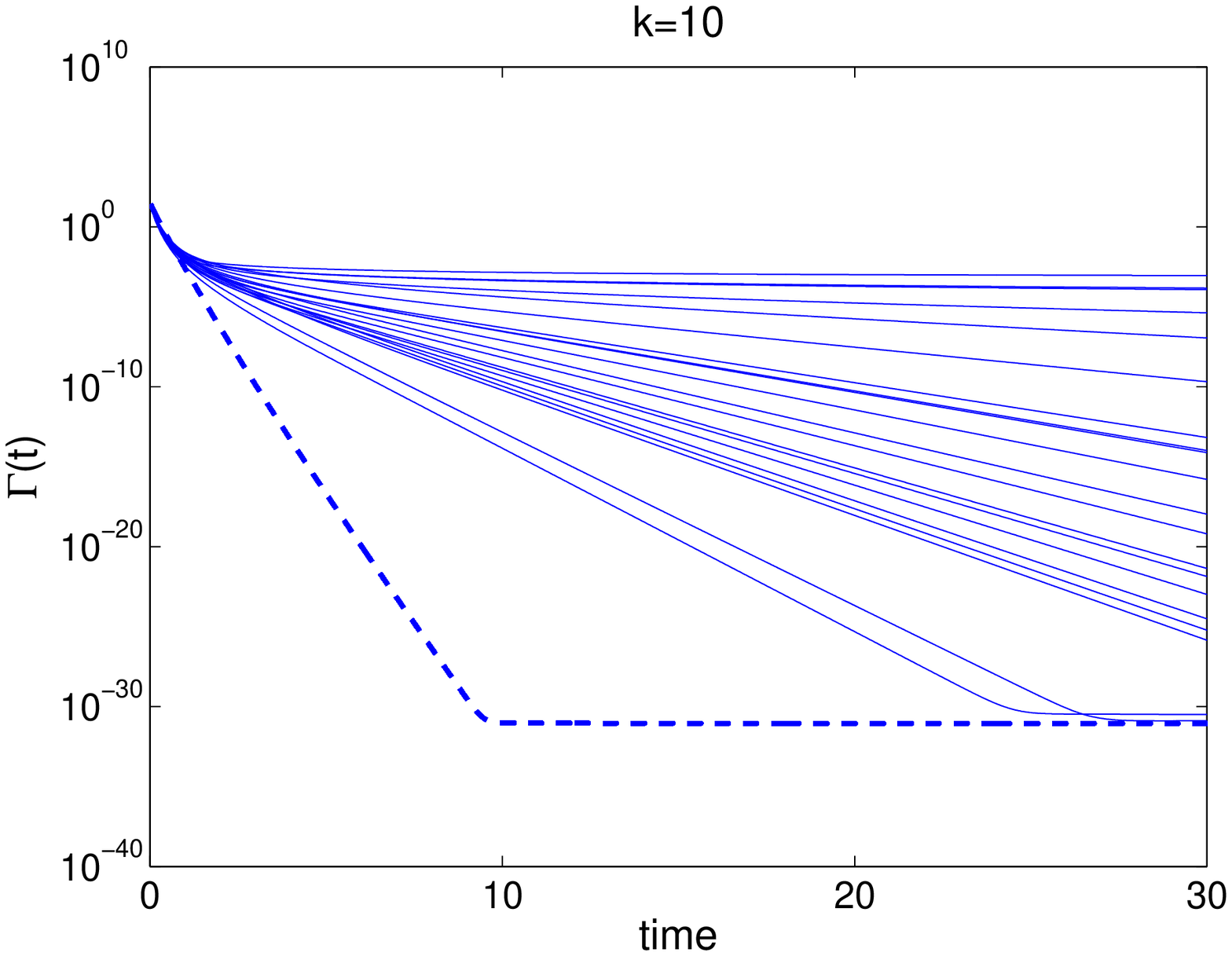}, width=0.5\linewidth} \\
      \epsfig{figure={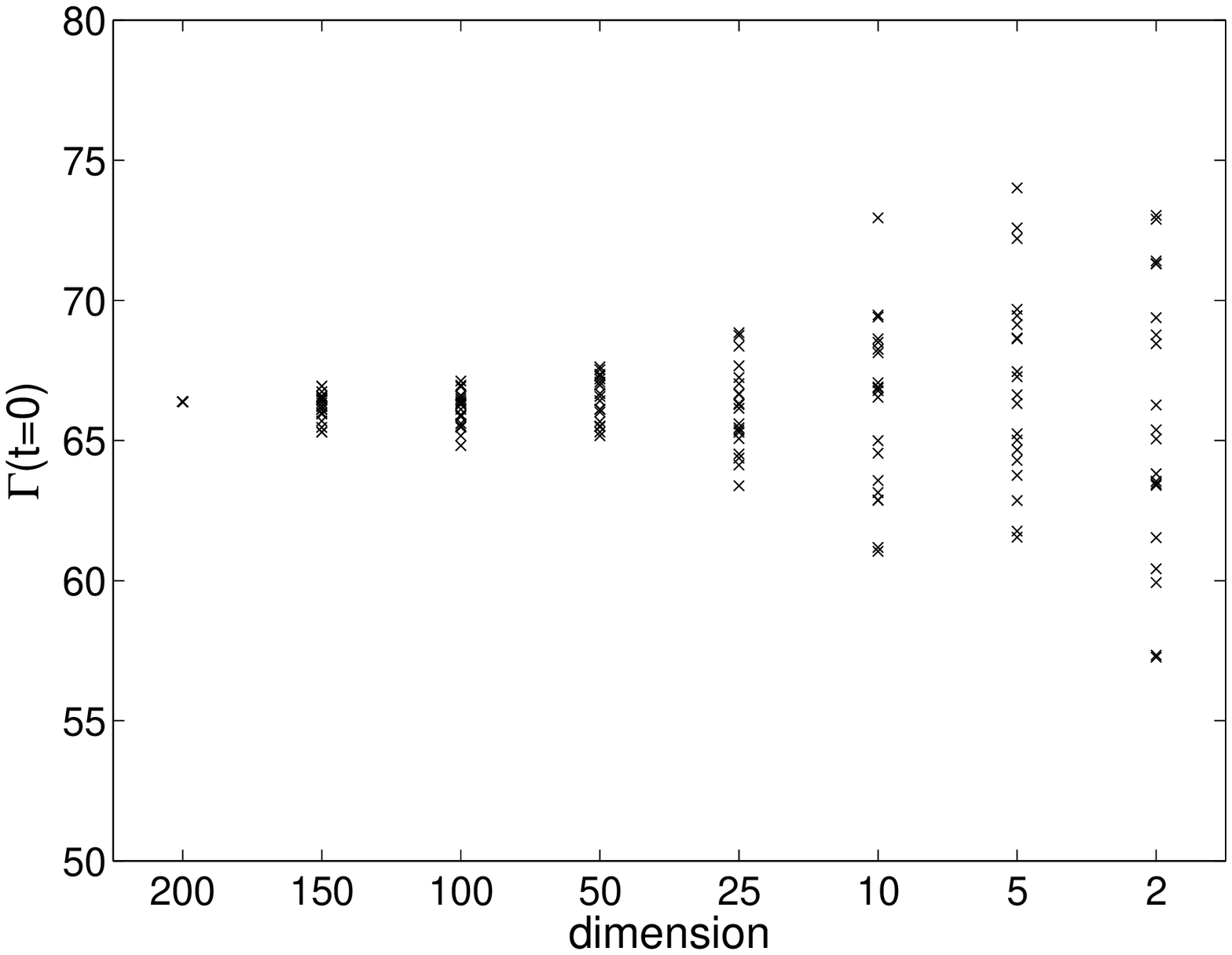}, width=0.5\linewidth} & &
      \epsfig{figure={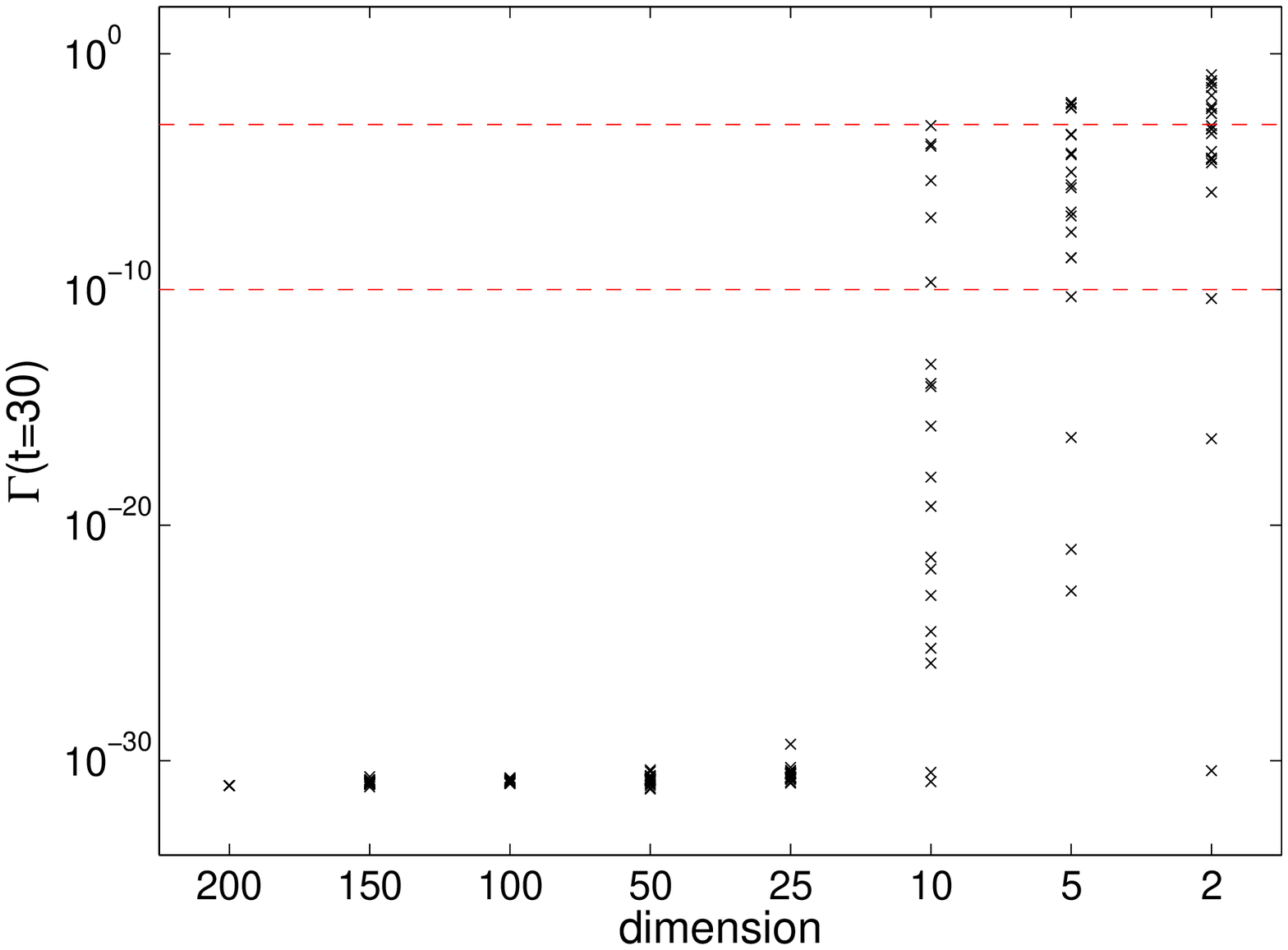}, width=0.5\linewidth}
\end{array}\]
\caption{Numerical results for $\beta=1.5$: First row shows the evolution of $\Gamma(t)$
of the system projected to dimension $k=100$ (left) and $k=10$ (right)
in the twenty realizations, compared to the original system (bold dashed line).
Second row shows the initial values $\Gamma(t=0)$ and final values $\Gamma(t=30)$
in all the performed simulations.}
\label{Numfig1}
\end{figure}

\begin{figure}[ht]
\[\begin{array}{ccc}
   \centering
      \epsfig{figure={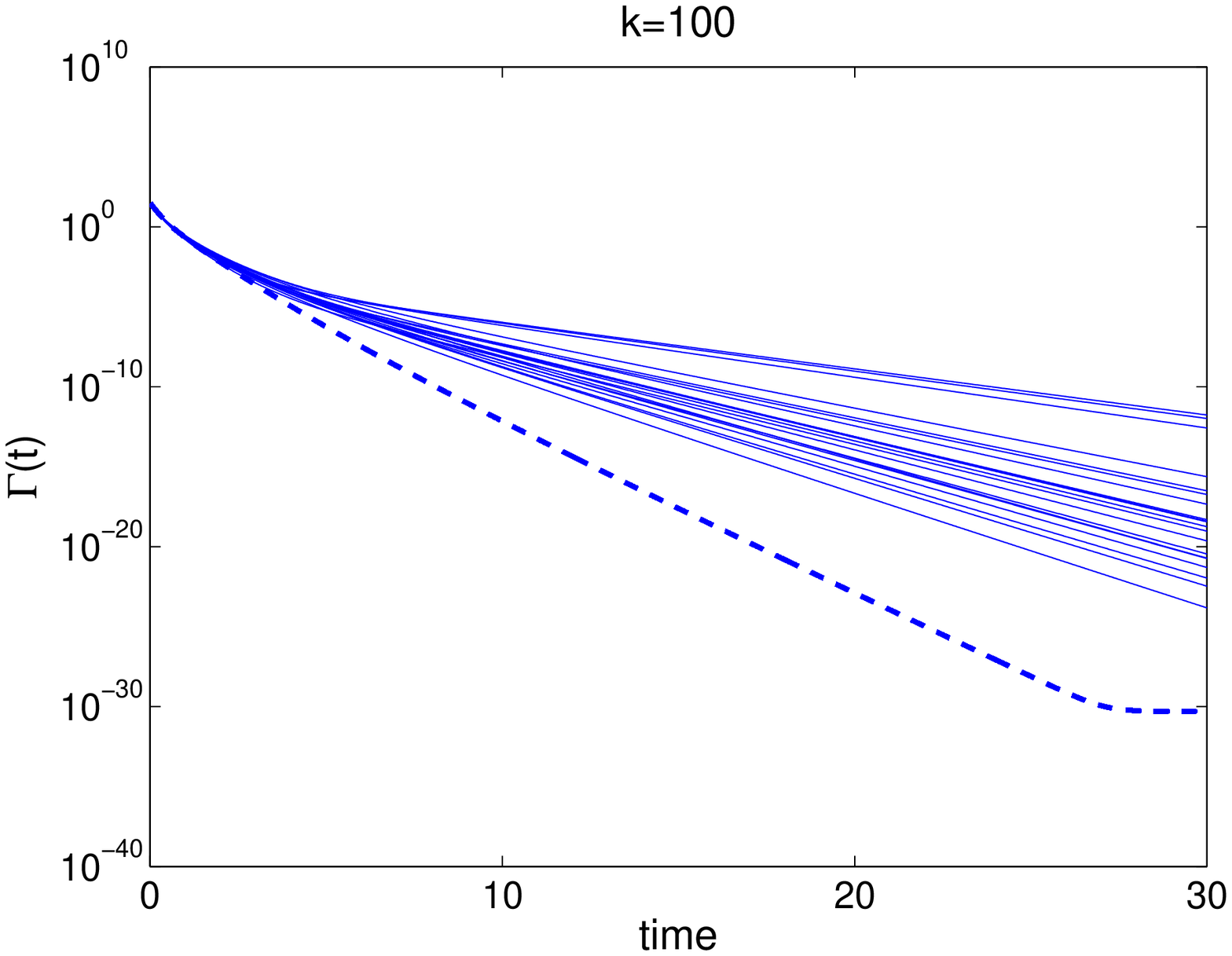}, width=0.5\linewidth} & &
      \epsfig{figure={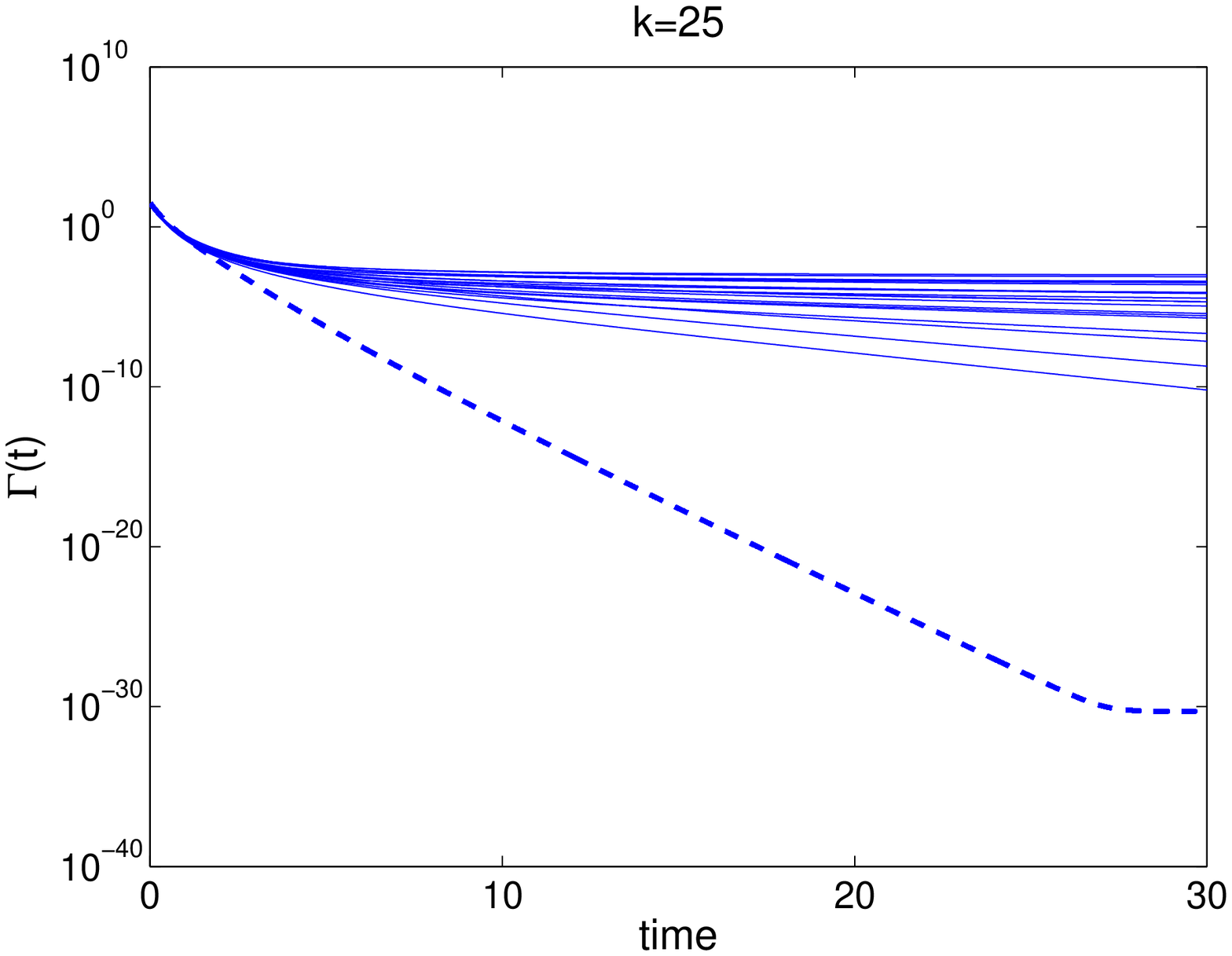}, width=0.5\linewidth} \\
      \epsfig{figure={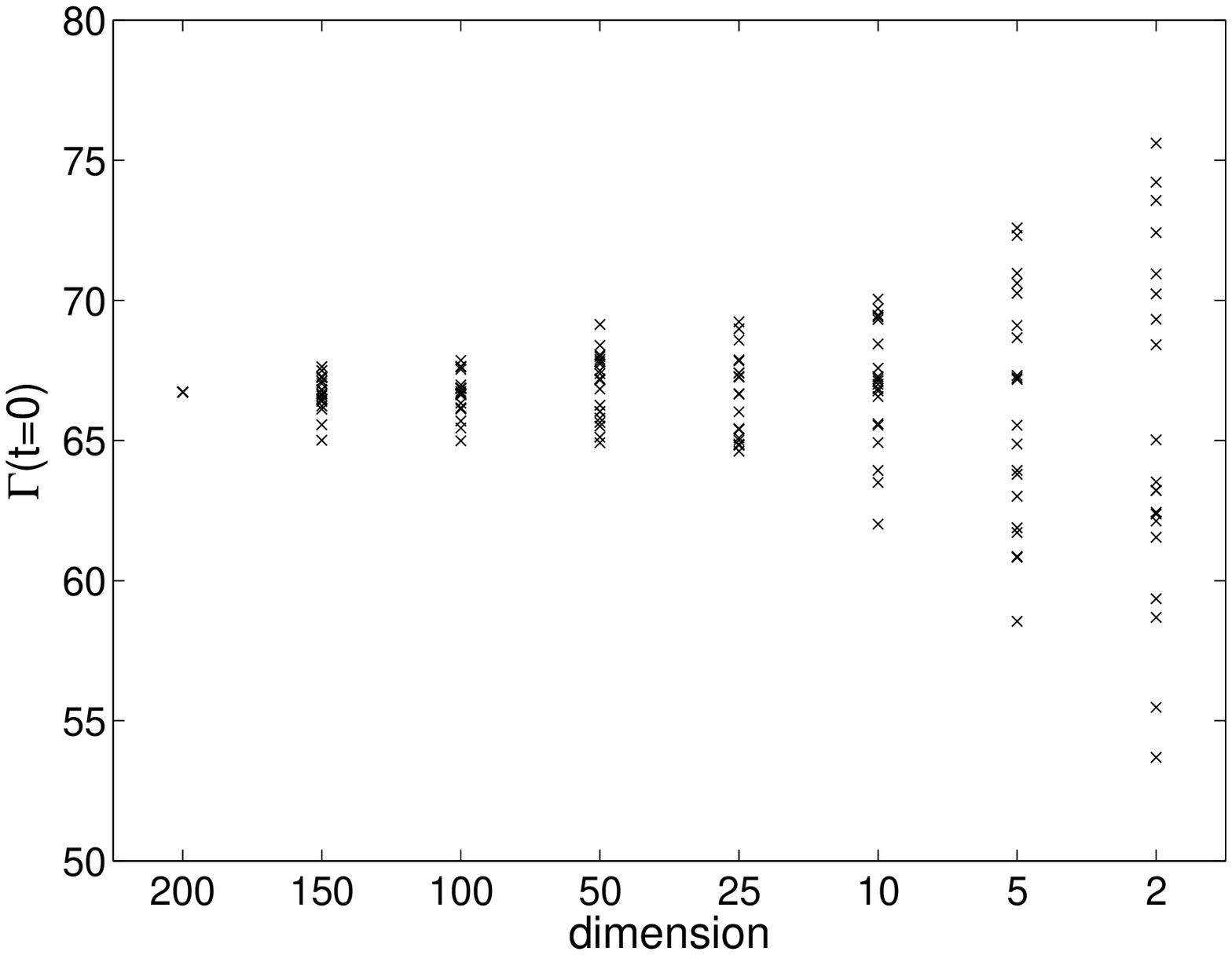}, width=0.5\linewidth} & &
      \epsfig{figure={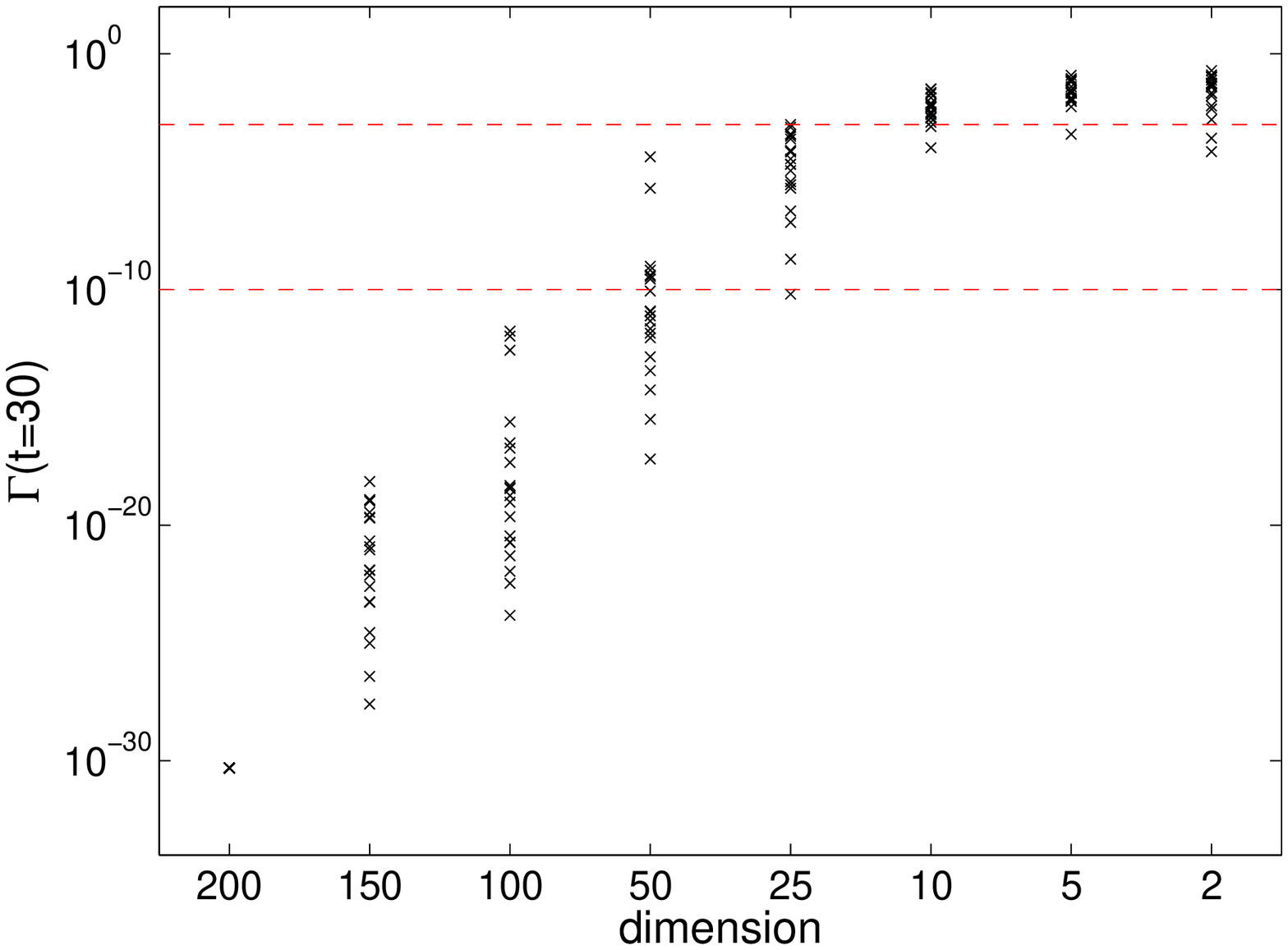}, width=0.5\linewidth}
\end{array}\]
\caption{Numerical results for $\beta=1.62$: First row shows the evolution of $\Gamma(t)$
of the system projected to dimension $k=100$ (left) and $k=25$ (right)
in the twenty realizations, compared to the original system (bold dashed line).
Second row shows the initial values $\Gamma(t=0)$ and final values $\Gamma(t=30)$
in all the performed simulations.}
\label{Numfig2}
\end{figure}

\begin{figure}[ht]
\[\begin{array}{ccc}
   \centering
      \epsfig{figure={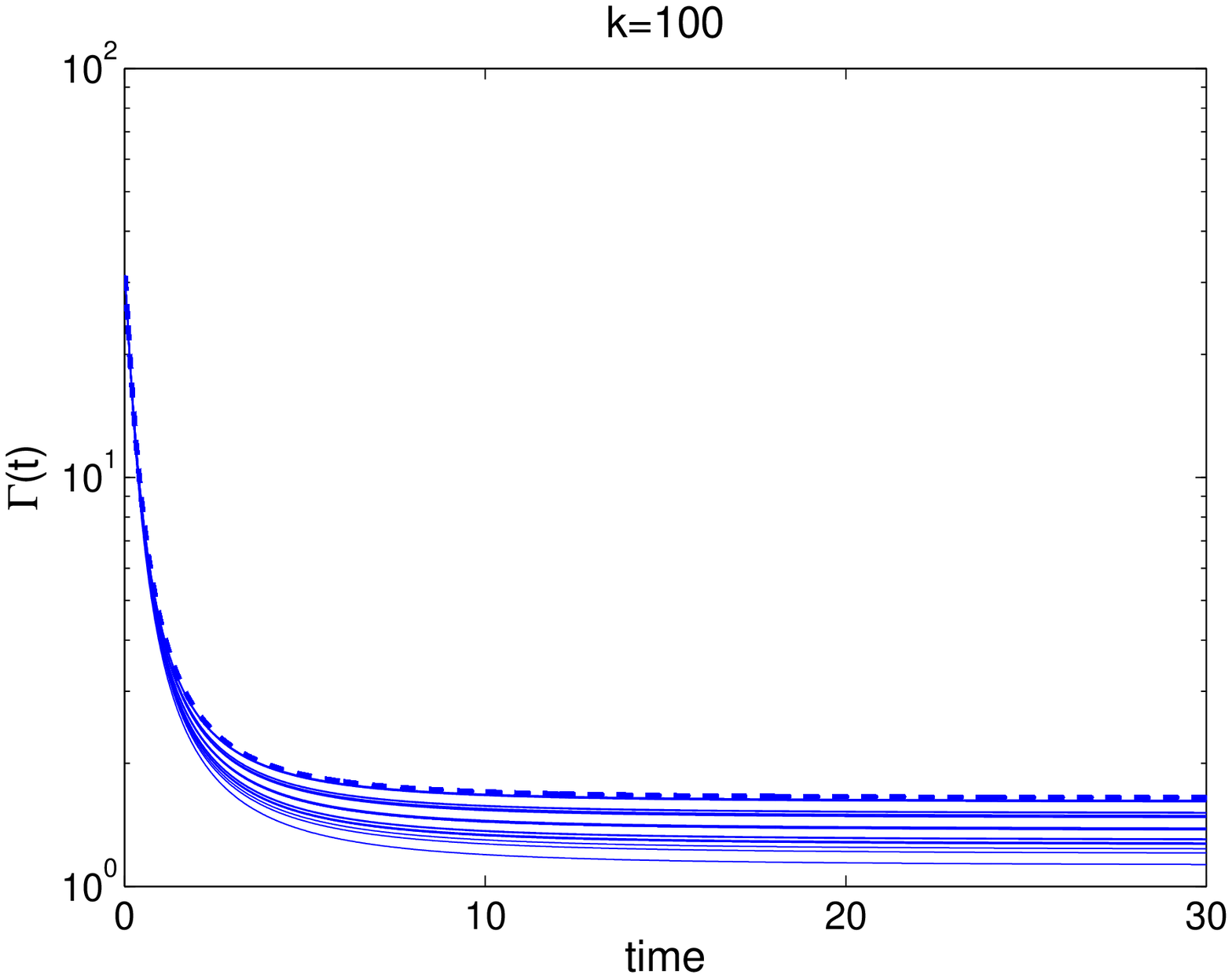}, width=0.5\linewidth} & &
      \epsfig{figure={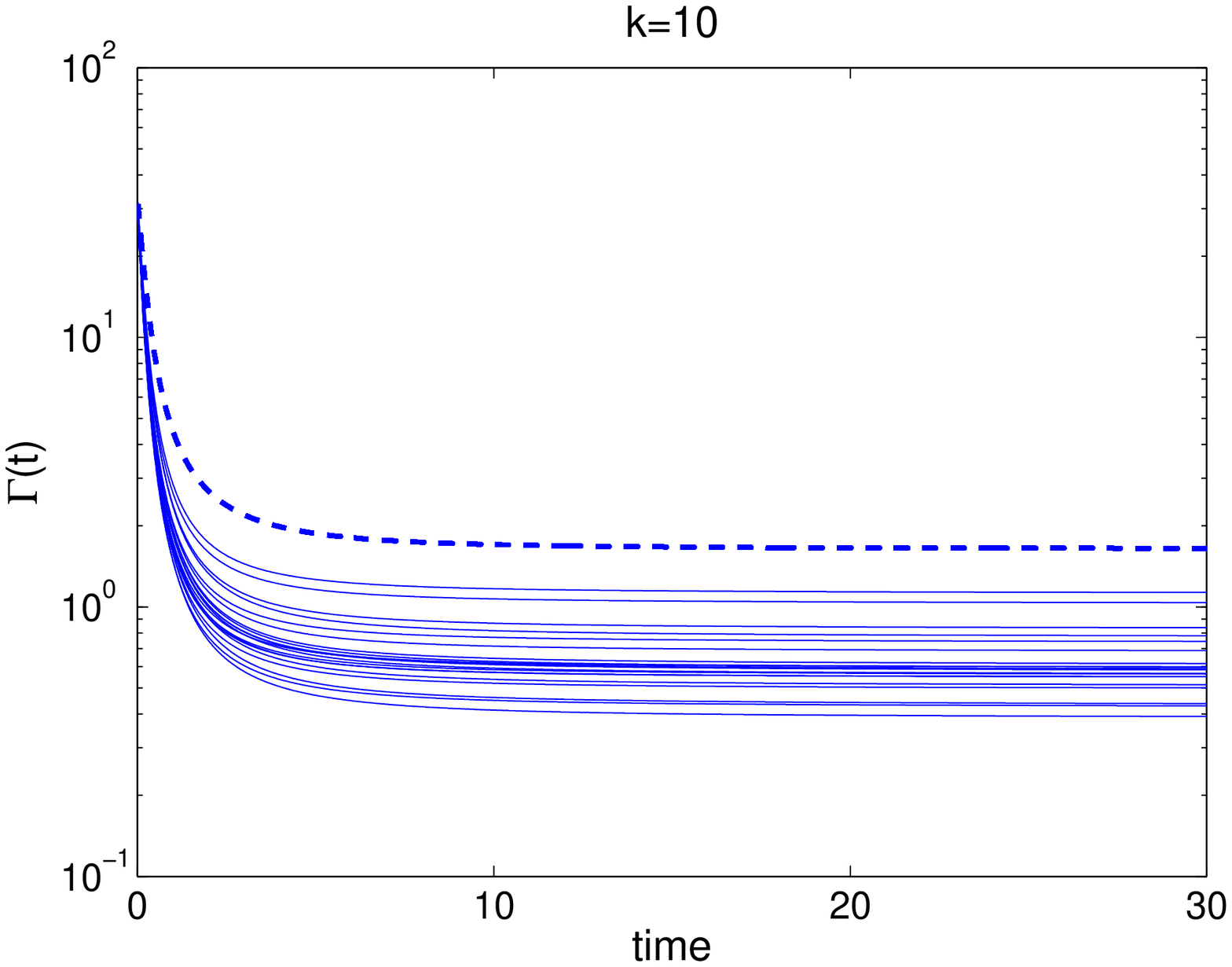}, width=0.5\linewidth} \\
      \epsfig{figure={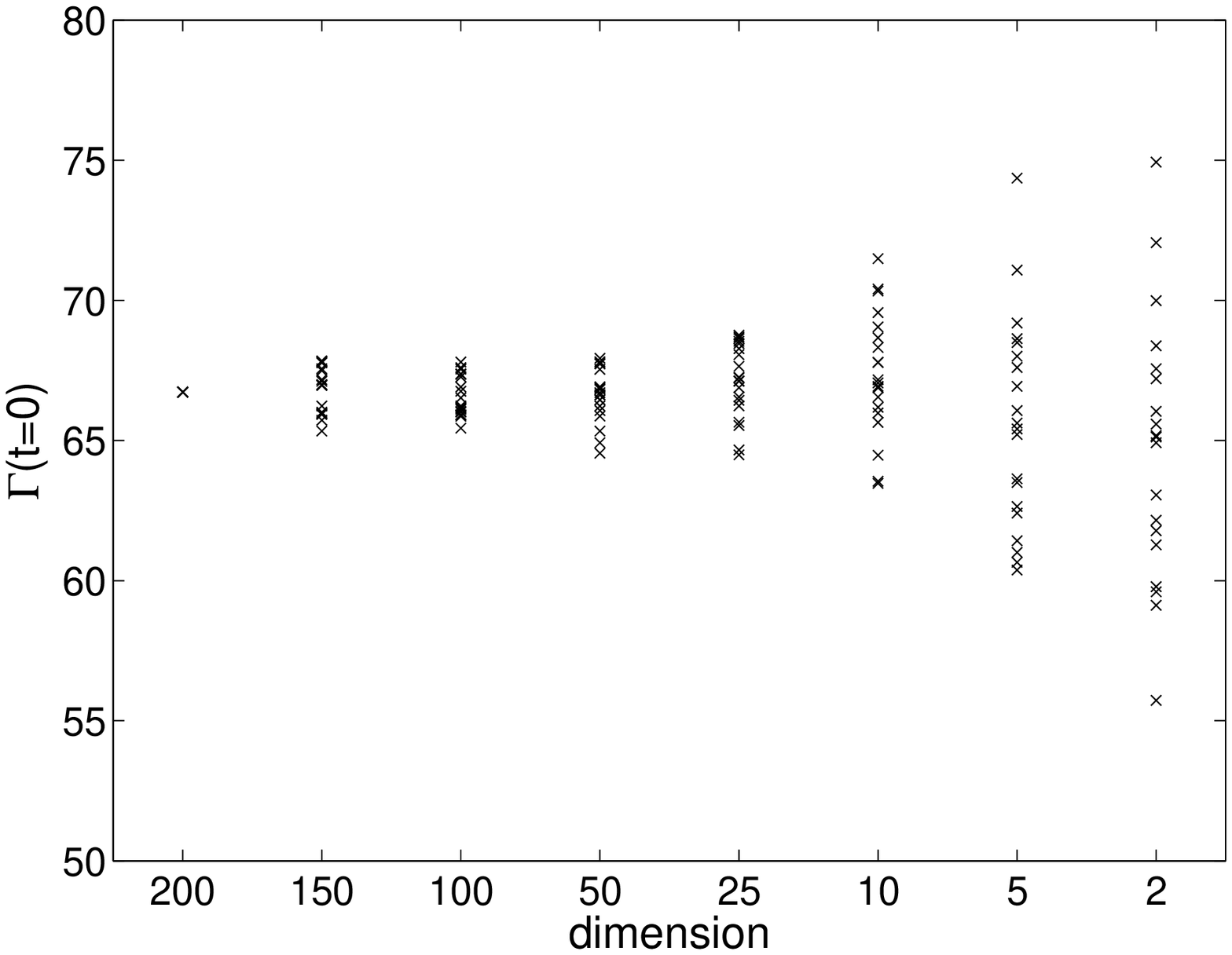}, width=0.5\linewidth} & &
      \epsfig{figure={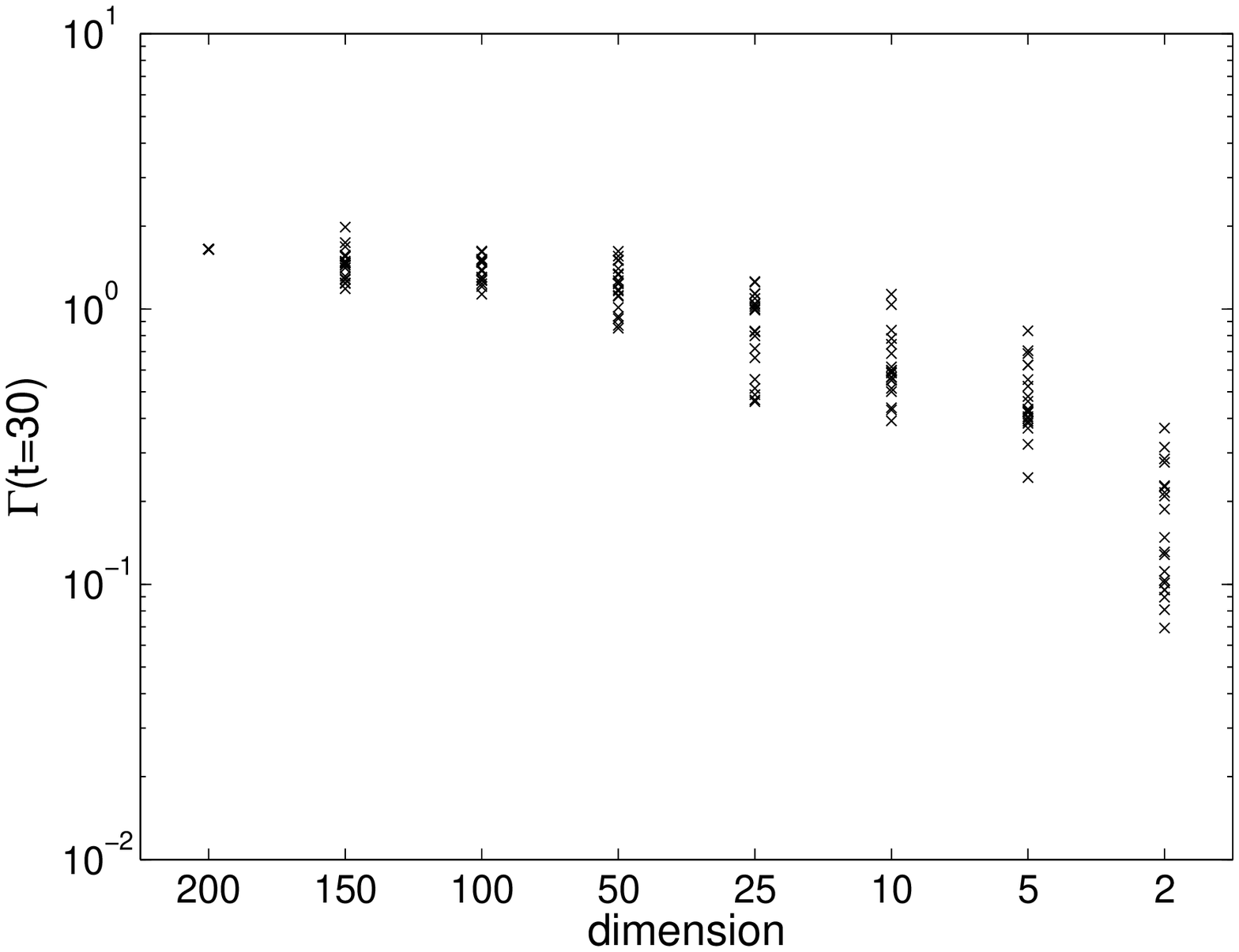}, width=0.5\linewidth}
\end{array}\]
\caption{Numerical results for $\beta=1.7$: First row shows the evolution of $\Gamma(t)$
of the system projected to dimension $k=100$ (left) and $k=10$ (right)
in the twenty realizations, compared to the original system (bold dashed line).
Second row shows the initial values $\Gamma(t=0)$ and final values $\Gamma(t=30)$
in all the performed simulations.}
\label{Numfig3}
\end{figure}

We illustrate this fact by a numerical experiment,
where we choose $N=1000$ and $d=200$,
i.e., every agent $i$ is determined by a $200$-dimensional vector $x_i$
of its state and a $200$-dimensional vector $v_i$
giving the rate of change of its state.
The initial datum $(x^0, v^0)$ is generated randomly,
every component of $x^0$ being drawn independently from the uniform distribution on $[0,1]$
and every component of $v^0$ being drawn independently from the uniform distribution on $[-1,1]$.
We choose $\beta = 1.5$, $1.62$ and $1.7$,
and with every of these values we perform the following set of simulations:
\begin{enumerate}
\item Simulation of the original system in $200$ dimensions.
\item Simulations in lower dimensions $k$: the initial condition $(x^0, v^0)$
is projected into the $k$-dimensional space by a random Johnson-Lindenstrauss
projection matrix $M$ with Gaussian entries. The dimension $k$
takes the values $150$, $100$, $50$, $25$, $10$, $5$, and $2$.
For every $k$, we perform the simulation twenty times,
each time with a new random projection matrix $M$.
\end{enumerate}
All the simulations were implemented in MATLAB, using $1500$ steps of the forward Euler method
with time step size $0.02$. The paths of $\Gamma(t)$ from the twenty experiments
with $k=100$ and $k=25$ or $k=10$ are shown in the first rows of Figs.~\ref{Numfig1},~\ref{Numfig2}
and, resp.,~\ref{Numfig3} for $\beta = 1.5$, $1.62$ and, resp., $1.7$.

The information we are actually interested in is whether flocking
takes place, in other words, whether the fluctuations of velocities $\Gamma(t)$
tend to zero.
Typically, after an initial phase, the graph of $\Gamma(t)$
gives a clear indication either about exponentially fast convergence to zero
(due to rounding errors, ``zero'' actually means values
of the order $~10^{-30}$ in the simulations)
or about convergence to a positive value.
However, in certain cases the decay may be very slow and
a very long simulation of the system would be needed to see if the limiting value
is actually zero or not. Therefore, we propose the following heuristic rules
to decide about flocking from numerical simulations:
\begin{itemize}
\item If the value of $\Gamma$ at the final time $t=30$
is smaller than $10^{-10}$, we conclude that flocking took place.
\item If the value of $\Gamma(30)$ is larger than $10^{-3}$,
we conclude that flocking did not take place.
\item Otherwise, we do not make any conclusion.
\end{itemize}
In the second rows of Figs.~\ref{Numfig1},~\ref{Numfig2} and~\ref{Numfig3} we present
the initial and final values of $\Gamma$ of the twenty simulations
for all the dimensions $k$, together with the original dimension $d=200$.
In accordance with the above rules, flocking takes place if the final
value of $\Gamma$ lies below the lower dashed line,
does not take place if it lies above the upper dashed line,
otherwise the situation is not conclusive.
The results are summarized in Table~\ref{Numtable1}.

\begin{table}[htb]
\footnotesize
\begin{tabular}{ccccc}
$\beta = 1.5$ & & $\beta = 1.62$ & & $\beta = 1.7$ \\
\begin{tabular}[h]{c|ccc}
dim & pos & neg & n/a \\
\hline
200 & 1 & 0 & 0 \\
150 & 20 & 0 & 0 \\
100 & 20 & 0 & 0 \\
50 & 20 & 0 & 0 \\
25 & 20 & 0 & 0 \\
10 & 14 & 0 & 6 \\
5 & 4 & 4 & 12 \\
2 & 3 & 8 & 9
\end{tabular}  & &
\begin{tabular}[h]{c|ccc}
dim & pos & neg & n/a \\
\hline
200 & 1 & 0 & 0 \\
150 & 20 & 0 & 0 \\
100 & 20 & 0 & 0 \\
50 & 13 & 0 & 7 \\
25 & 1 & 1 & 18 \\
10 & 0 & 18 & 2 \\
5 & 0 & 19 & 1 \\
2 & 0 & 18 & 2
\end{tabular}  & &
\begin{tabular}[h]{c|ccc}
dim & pos & neg & n/a \\
\hline
200 & 0 & 1 & 0 \\
150 & 0 & 20 & 0 \\
100 & 0 & 20 & 0 \\
50 & 0 & 20 & 0 \\
25 & 0 & 20 & 0 \\
10 & 0 & 20 & 0 \\
5 & 0 & 20 & 0 \\
2 & 0 & 20 & 0
\end{tabular}
\end{tabular}
\caption{Statistics of the flocking behaviors of the systems in the original dimension $d=200$ and
in the projected dimensions. With $\beta=1.5$ and $\beta=1.62$, the original system ($d=200$)
exhibited flocking behavior. With $\beta=1.5$, even after random projections into $25$ dimensions,
the system exhibited flocking in all $20$ repetitions of the experiment, and still in $14$ cases in dimension $10$.
With $\beta=1.62$, the deterioration of the flocking behavior with decreasing dimension was much faster,
and already in dimension $25$ the situation was not conclusive. This is related to the fact that the value
$\beta=1.62$ was chosen to intentionally bring the system close to the borderline
between flocking and non-flocking. Finally, with $\beta=1.7$, the original system did not flock,
and, remarkably, all the projected systems (even to two dimensions) exhibit the same behavior.
}
\label{Numtable1}
\end{table}

Experience gained with a large amount of numerical experiments
shows the following interesting fact: The flocking behavior
of the Cucker-Smale system is very stable with respect
to the Johnson-Lindenstrauss projections.
Usually, the projected systems show the same flocking behavior
as the original one, even if the dimension is reduced dramatically,
for instance from $d=200$ to $k=10$
(see Figs~\ref{Numfig1} and~\ref{Numfig3}).
This stability can be roughly explained as follows:
Since the flocking behavior depends mainly on the initial
values of $\Gamma$ and $\Lambda$, which are statistical properties
of the random distributions used for the generation of initial data,
and since $N$ is sufficiently large, the concentration
of measure phenomenon takes place. Its effect is
that the initial values of the fluctuations of the projected data
are very close to the original ones, and thus the flocking behavior is
(typically) the same.
There is only a narrow interval of values of $\beta$ (in our case this
interval is located around the value $\beta= 1.62$),
which is a borderline region between flocking and non-flocking,
and the projections to lower dimensions spoil the flocking behavior,
see Fig~\ref{Numfig2}.
Let us note that in our simulations we were only able
to detect cases when flocking took place in the original system,
but did not take place in some of the projected ones.
Interestingly, we never observed the inverse situation,
a fact which we are not able to explain satisfactorily.
In fact, one can make other interesting observations, deserving
further investigation. For instance, Figs.~\ref{Numfig1} and~\ref{Numfig2}
show that if the original system exhibits flocking, then the
curves of $\Gamma(t)$ of the projected systems tend to lie above the
curve of $\Gamma(t)$ of the original one. The situation is reversed
if the original system does not flock, see Fig.~\ref{Numfig3}.

From a practical point of view, we can make the following conclusion:
To obtain an indication about the flocking behavior of a highly dimensional Cucker-Smale system,
it is typically satisfactory to perform a limited number of simulations
of the system projected into a much lower dimension, and evaluate the statistics
of their flocking behavior. If the result is the same for the majority
of simulations, one can conclude that the original system
very likely has the same flocking behavior as well.

\begin{figure}[ht]
\[\begin{array}{ccc}
   \centering
      \epsfig{figure={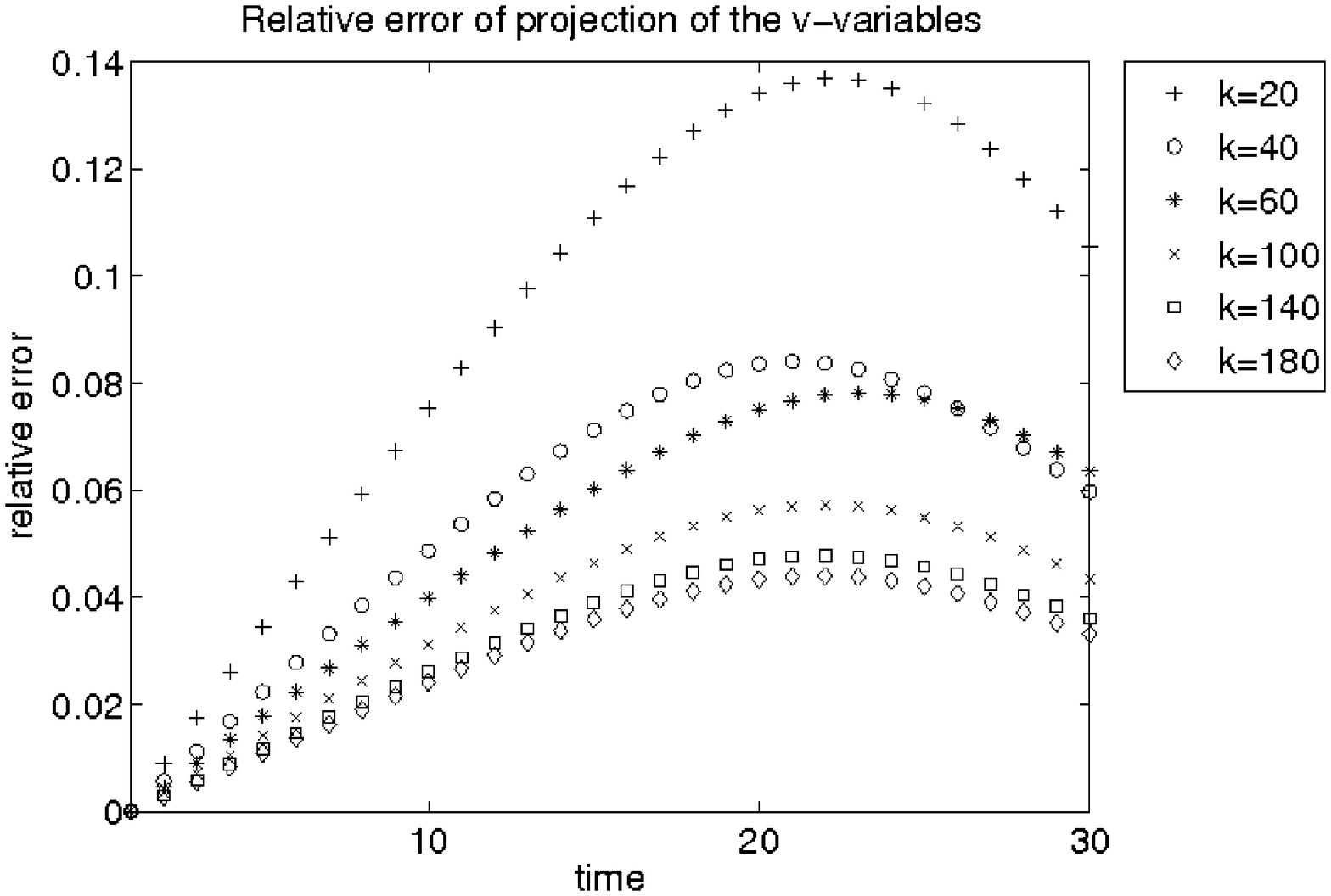}, width=0.52\linewidth} & &
      \epsfig{figure={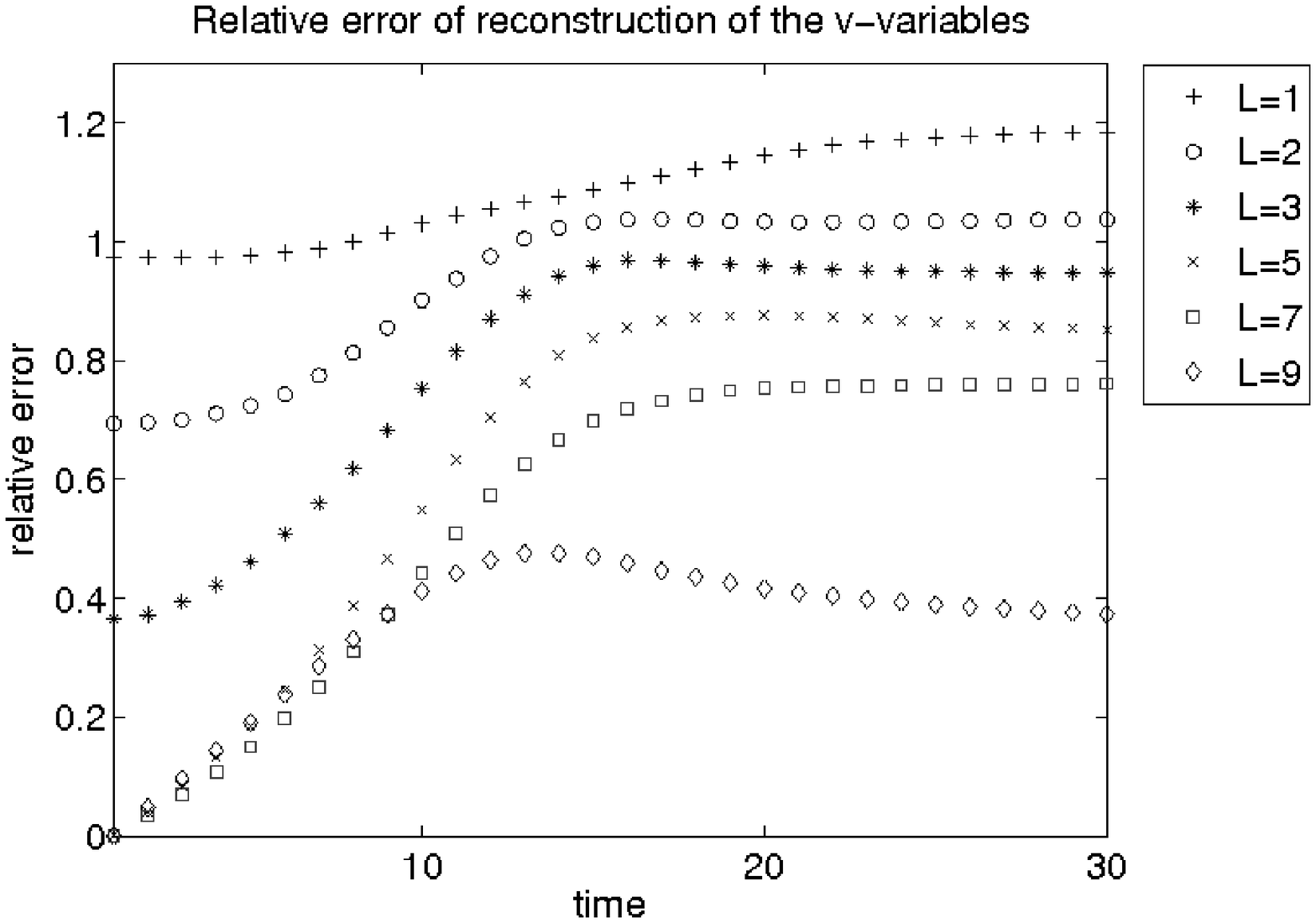}, width=0.5\linewidth}
\end{array}\]
\caption{Numerical results showing the time evolution of the relative error of projection (left panel)
and relative error of recovery via $\ell_1$-minimization (right panel) of the $v$-variables.
}
\label{Numfig4}
\end{figure}
\subsubsection{Numerical validation of the high and low dimensional approximation properties}
Finally, we show how the relative error of projection and recovery evolves in time.
We consider an initial datum $(x^0,v^0)\in\R^{N\times d}\times\R^{N\times d}$ for the Cucker-Smale system
with $N=d=200$ and randomly generated entries from the normal distribution.
The parameter $\beta=0.4$, therefore, the system will exhibit flocking.
First we project the system into $k=20,40,60,100,140,180$ dimensions
and calculate the relative error of the projection of the $v$-variables, given by
\[
    \left(\frac{\sum_{i=1}^N \|Mv_i-v_j\|_{\ell_2^k}^2}{\sum_{i=1}^N \|Mv_i\|_{\ell_2^k}^2} \right)^{1/2} \,.
\]
We observe that the maximal relative error (for $k=20$) is around $14\%$, which we consider as a very
good result.
Moreover, in all $9$ cases, the error first increases, but after $t\simeq 22$ it starts decreasing,
which is a consequence of the flocking behavior and concentration of measure, see the graphics in 
Figure~\ref{Numfig4} on the left.
This clearly shows that the worst-case estimate of Theorem~\ref{thm1} with exponential growth in time
is overly pessimistic.

In our second experiment, we take a randomly generated initial condition with $N=d=200$
and $80\%$ of the entries set to zero.
Then, we take $L$ projections of the system into $20$ dimensions, with $L=1,2,3,5,7,9$,
and reconstruct the $v$-trajectories using $\ell_1$-minimization, as described
in Section~\ref{CSsec}.
In the graphics in Figure~\ref{Numfig4} on the right, we plot the relative errors, given by
\[
    \left(\frac{\sum_{i=1}^N \|\tilde v_i-v_i\|_{\ell_2^k}^2}{\sum_{i=1}^N \|v_i\|_{\ell_2^k}^2} \right)^{1/2} \,,
\]
where $\tilde v_i$ are the recovered trajectories.
Again, we observe that the errors grow much slower than exponentially,
and after $t\simeq 15$ they even tend to stay constant or slightly decrease.

\section{Mean-field limit and kinetic equations in high dimension}
In the previous sections we were concerned with tractable simulation of the dynamical systems of the type~\eqref{eq:dyn}
when the dimension $d$ of the parameter space is large.
Another source of possible intractability in numerical simulations appears in the situation where the number of agents $N$ is very large.
In general, large $N$ imposes even a much more severe limitation than large $d$, since the computational complexity of~\eqref{eq:dyn} is $\O(d\times N^2)$.
Therefore, in the next sections we consider the so-called {\it mean-field limit} of~\eqref{eq:dyn} as $N\to\infty$,
where the evolution of the system is described by time-dependent probability measures $\mu(t)$ on $\mathbb R^d$, representing the density distribution of agents, and satisfying mesoscopic partial differential equations of the type \eqref{genkin}.
This strategy originated from the kinetic theory of gases, see \cite{cip} for classical references.
We show how our projection method can be applied for dimensionality reduction of the corresponding kinetic equations
and explain how the probability measures can be approximated by atomic measures.
Using the concepts of {\it delayed curse of dimension} and {\it measure quantization} known from optimal integration problems in high dimension,
we show that under the assumption that the measure concentrates along low-dimensional subspaces (and more generally along low-dimensional sets or manifolds),
it can be approximated by atomic measures with sub-exponential (with respect to $d$) number of atoms.
Through such approximation, we shall show that we can approximate suitable random averages of the solution of the original partial differential equation in high dimension by tractable simulations
of corresponding solutions of lower-dimensional kinetic equations.

Another interesting approach to the problem of efficient numerical simulation of large group dynamics
is the so-called ``equation-free'' approach, see e.g.~\cite{mnllk}.
Here, convenient coarse-grained variables that account
for rapidly developing correlations during initial transients are chosen,
in order to perform efficient computations of coarse-grained steady states
and their bifurcation analysis.
The big advantage of the equation-free approach is that the coarse-grained dynamics can be explored
without the assumption of the continuum limit equation as we consider here.
The premise of the method is that coarse-grained
governing equations conceptually exist, but are
not explicitly available in closed form. 
The main idea is that short bursts of appropriately
initialized microscopic (fine-scale) simulations
and the projection of the results onto coarse-grained variables
result in time-steppers (mappings) for those variables (which is
effectively the same as the discretization of the unavailable
equations).
One then processes the results of the short simulations to
estimate various coarse-grained quantities (such as time
derivatives, action of Jacobians, residuals) to perform relevant coarse-grained
level numerical computations, as if
those quantities were obtained from coarse-grained governing equations.

\subsection{Formal derivation of mean-field equations}
In this section we briefly explain how the mean-field limit description corresponding to~\eqref{eq:dyn}
can be derived.
This is given, under suitable assumptions on the family of the governing functions
$\mathcal F_N =\{f_i, f_{ij}: i,j =1, \dots N\}$, by the general formula
\begin{equation}\label{genkin}
\frac{\partial \mu}{\partial t} + \nabla \cdot ( \mathcal H_{\mathcal F}[\mu] \mu) =0,
\end{equation}
where $\mathcal H_{\mathcal F}[\mu]$ is a field in $\R^d$,
determined by the sequence $\mathcal F=(\mathcal F_N)_{N \in \mathbb N}$.

In order to provide an explicit example,
we show how to formally derive the mean field limit of systems of the type
\begin{eqnarray}
    \dot x_i &=& v_i \,, \label{eq:dyn4}\\
    \dot v_i &=& \sum_{j=1}^N f^{vv}_{ij}(\D x, \D v) v_j + \sum_{j=1}^N f^{vx}_{ij}(\D x) x_j\,,\label{eq:dyn5}
\end{eqnarray}
with
\[
   f^{vx}_{ij}(\D x) &=& - \frac{\delta_{ij}} {N}\sum_{k\neq i} u(\|x_i-x_k\|_{\ell_2^d})
           + \frac{1- \delta_{ij}} {N} u(\|x_i-x_j\|_{\ell_2^d}) \,,\\
   f^{vv}_{ij}(\D x,\D v) &=& \delta_{ij} \left ( h(\|v_i\|_{\ell_2^d}^2)
        - \frac{1} {N} \sum_{k=1}^N g(\|x_i-x_k\|_{\ell_2^d}) \right)
        + \frac{1- \delta_{ij}} {N} g(\|x_i-x_j\|_{\ell_2^d}) \,.
\]
Note that for suitable choices of the functions $h,g,u$
this formalism includes both the Cucker-Smale model~\eqref{CS1}--\eqref{CS2}
and the self-propulsion and pairwise interaction model~\eqref{DOrsogna1}--\eqref{DOrsogna2}.
We define the empirical measure associated to the solutions $x_i(t)$, $v_i(t)$
of \eqref{eq:dyn4}--\eqref{eq:dyn5} as
$$
   \mu^N(t):=\mu^N(t,x,v)= \frac{1}{N} \sum_{i=1}^N \delta_{x_i(t)}(x) \delta_{v_i(t)}(v) \,.
$$
Taking a smooth, compactly supported test function $\xi\in C^\infty_0(\R^{2d})$
and using \eqref{eq:dyn4}--\eqref{eq:dyn5},
one easily obtains by a standard formal calculation (see~\cite{CFTV})
\begin{eqnarray}
   \frac{d}{dt} \langle \mu^N(t) , \xi \rangle 
  &=& \frac{d}{dt} \left (\frac{1}{N} \sum_{i=1}^N \xi(x_i(t),v_i(t)) \right) \label{BBGKY}\\
  &=& \int_{\R^{2d}} \grad_x\xi(x,v)\cdot v \d \mu^N(t,x,v)
     + \int_{\R^{2d}} \grad_v\xi(x,v)\cdot \mathcal{H}[\mu^N(t)](x,v) \d \mu^N(t,x,v) \nonumber \,,
\end{eqnarray}
with
\[
   \mathcal{H}[\mu](x,v) = h(\|v\|_{\ell_2^d}) v
       + \int_{\R^{2d}} g(\|x-y\|_{\ell_2^d}) (w-v) \d\mu(y,w)
       + \int_{\R^{2d}} u(\|x-y\|_{\ell_2^d}) (y-x) \d\mu(y,w) \,.
\]
We now assume weak convergence of a subsequence of $(\mu^N(t))_{N \in \mathbb N}$
to a time-dependent measure $\mu(t)=\mu(t,x,v)$ and boundedness of its first order moment,
which indeed can be established rigorously for the Cucker-Smale and the self-propulsion and pairwise interaction systems
(see~\cite{TH},~\cite{DOrsogna}).
Then, passing to the limit $N\to\infty$ in~\eqref{BBGKY}, one obtains in the strong formulation
that $\mu$ is governed by
$$
   \frac{\partial \mu}{\partial t}(t,x,v) +  v \cdot \nabla_x \mu(t,x,v)
          + \nabla_v \cdot \left ( \mathcal H[\mu(t)](x,v) \mu(t,x,v) \right ) = 0 \,,
$$
which is an instance of the general prototype \eqref{genkin}.

Using the same formal arguments as described above,
one can easily derive mean field limit equations corresponding to~\eqref{eq:dyn}
with different choices of the family $\mathcal F$.

\subsection{Monge-Kantorovich-Rubinstein distance and stability}
In several relevant cases, and specifically for the Cucker-Smale and the self-propulsion and pairwise interaction systems \cite{CCR},
solutions of equations of the type \eqref{genkin}
are stable with respect to suitable distances.
We consider the space $\mathcal P_1(\mathbb R^d)$,
consisting of all probability measures on $\mathbb R^d$ with finite first moment.
In $\mathcal P_1(\mathbb R^d)$ and for solutions of \eqref{genkin},
a natural metric to work with is the so-called {\it Monge-Kantorovich-Rubinstein distance} (also called {\it Wasserstein distance}) ~\cite{Villani},
\begin{equation}\label{wasserstein}
W_1(\mu,\nu) := \sup \{ \left | \langle \mu - \nu, \xi \rangle \right |= \left| \int_{\mathbb R^d} \xi(x) d(\mu-\nu)(x) \right |, \xi \in {\rm Lip}(\mathbb R^d),  {\rm Lip}(\xi) \leq 1 \}.
\end{equation}
We further denote $\mathcal P_c(\mathbb R^d)$ the space of compactly supported probability measures on $\mathbb R^d$.
In particular, throughout the rest of this paper,
we will assume that for any compactly supported measure valued weak solutions
$\mu(t),\nu(t) \in C([0,T],\mathcal P_c(\mathbb R^d))$ of \eqref{genkin}
we have the following stability inequality
\begin{equation}\label{stability}
W_1(\mu(t), \nu(t)) \leq C(t) W_1(\mu(0), \nu(0)), \quad t \in [0,T],
\end{equation}
where $C(t)$ is a positive increasing function of $t$ with $C(0)>0$, independent of the dimension $d$.
We address the interested reader to \cite[Section 4]{CCR} for a sample of general conditions
on the vector field $\mathcal H[\mathcal F](\mu)$ which guarantee stability \eqref{stability}
for solutions of equations \eqref{genkin}.

\subsection{Dimensionality reduction of kinetic equations}
Provided a high-dimensional measure valued solution to the equation
\begin{equation}\label{genkin2}
\frac{\partial \mu}{\partial t} + \nabla \cdot ( \mathcal H_{\mathcal F}[\mu] \mu) =0, \quad \mu(0)= \mu_0 \in \mathcal P_c(\mathbb R^d) \,,
\end{equation}
we will study the question whether its solution can be approximated by suitable projections in lower dimension.

Given a probability measure $\mu \in \mathcal P_1(\mathbb R^d)$,
its projection into $\R^k$ by means of a matrix $M:\mathbb R^d \to \mathbb R^k$
is given by the {\it push-forward} measure $\mu_M := M \# \mu$,
\begin{equation}\label{pf}
\langle \mu_M, \varphi \rangle := \langle \mu, \varphi(M \cdot) \rangle \quad \mbox{for all } \varphi \in {\rm Lip}(\mathbb R^{k}).
\end{equation}
Let us mention two explicit and relevant examples:
\begin{itemize}
\item If $\mu^N= \frac{1}{N}\sum_{i=1}^N \delta_{x_i}$ is an atomic measure,
we have $\langle \mu_M^N, \varphi \rangle = \langle \mu^N, \varphi(M \cdot) \rangle = \frac{1}{N}\sum_{i=1}^N \varphi(M x_i)$.
Therefore,
\begin{equation}\label{projatom}
\mu_M^N = \frac{1}{N}\sum_{i=1}^N \delta_{M x_i}\,.
\end{equation}
\item If $\mu$ is absolutely continuous with respect to the Lebesgue measure, i.e., it is a function in $L^1(\mathbb R^d)$, the calculation requires a bit more effort:
Let us consider $M^\dagger$ the pseudo-inverse matrix of $M$. Recall that $M^\dagger = M^* (M M^*)^{-1}$ is a right inverse of $M$, and $M^\dagger M$ is the orthogonal projection onto the range of $M^*$.
Moreover, $x= M^\dagger M x + \xi_x$, where $\xi_x \in \ker M$ for all $x \in \mathbb R^d$.
According to these observations, we write
\begin{eqnarray*}
\int_{\mathbb R^d} \varphi (M x) \mu (x)dx &=&  \int_{\mathbb R^d} \varphi (M x) \mu (M^\dagger M x + \xi_x)dx\\
&=& \int_{{\rm ran} M^* \oplus \ker M} \varphi (M x) \mu (M^\dagger M x + \xi_x)dx\\
&=& \int_{{\rm ran} M^* } \int_{\ker M} \varphi (M v) \mu (M^\dagger M v + v^\perp)dv^\perp dv
\end{eqnarray*}
Note now that $M_{| {\rm ran} M^*} : {\rm ran} M^*  \to {\rm ran} M \eqsim \mathbb R^k$ is an isomorphism,
hence $y = M v$ implies the change of variables $dv= \det(M_{| {\rm ran} M^*})^{-1} d y
= \det(M M^*)^{-1/2} d y$.
Consequently, we have
\begin{eqnarray*}
\int_{\mathbb R^d} \varphi (M x) \mu (x)dx &=&  \int_{\mathbb R^d} \varphi (M x) \mu (M^\dagger M x + \xi_x)dx\\
&=& \int_{{\rm ran} M^* } \int_{\ker M} \varphi (M v) \mu (M^\dagger M v + v^\perp)dv^\perp dv\\
&=& \int_{\mathbb R^k} \left (\frac{1}{\det(M M^*)^{1/2}} \int_{\ker M}  \mu (M^\dagger y + v^\perp)dv^\perp \right) \varphi (y) dy \,,
\end{eqnarray*}
and
$$
\mu_M(y) = \frac{1}{\det(M M^*)^{1/2}} \int_{\ker M}  \mu (M^\dagger y + v^\perp)dv^\perp.
$$
\end{itemize}
According to the notion of push-forward, we can consider the measure valued function
$\nu \in C([0,T],\mathcal P_c (\mathbb R^k))$, solution of the equation
\begin{equation}\label{genkin3}
\frac{\partial \nu}{\partial t} + \nabla \cdot ( \mathcal H_{\mathcal F_M}[\nu] \nu) =0, \quad \nu(0)= (\mu_0)_M \in \mathcal P_c(\mathbb R^k),
\end{equation}
where $(\mu_0)_M =  M \# \mu_0$ and $\mathcal F_M= (\{M f_i, f_{ij}, i,j=1,\dots, N\})_{N \in \mathbb N}$. As for the dynamical system \eqref{lowdimsys}, also equation \eqref{genkin3} is fully defined on the lower-dimensional space $\mathbb R^k$ and depends on the original high-dimensional problem exclusively by means of the initial condition.

The natural question at this point is whether the solution $\nu$ of \eqref{genkin3}
provides information about the solution $\mu$ of \eqref{genkin2}.
In particular, similarly to the result of Theorem~\ref{Thm:Cont},
we will examine whether the approximation 
$$
\nu(t) \approx \mu_M(t), \quad t \in [0,T],
$$
in Monge-Kantorovich-Rubinstein distance is preserved in finite time. We depict the expected result by the following diagram:
\begin{equation*}
\begin{matrix}
&\mu(0) & \stackrel{t}{\longrightarrow} & \mu(t) & \cr &\downarrow M & &\downarrow M& \cr &
\nu(0)=(\mu_0)_M&\stackrel{t}{\longrightarrow} &\nu(t) \approx \mu_M(t)&.  
\end{matrix}
\end{equation*}
This question will be addressed by approximation of the problem
by atomic measures and by an application of Theorem~\ref{Thm:Cont} for the corresponding dynamical system, as concisely described by
\begin{equation*}
\begin{matrix}
&\mu & \stackrel{W_1(\mu,\,\mu^N) \lesssim \varepsilon}{\longrightarrow} & \mu^N & \cr & \downarrow M & &  \downarrow M& \cr &
\nu \approx \mu_M &\stackrel{W_1(\nu,\,\nu^N) \lesssim \varepsilon}{\longrightarrow} &  \nu^N \approx \mu^N_M &  
\end{matrix}
\end{equation*}
Let us now recall the framework and general assumptions for this analysis to be performed.
We assume again that for all $N\in \mathbb N$ the family $\mathcal F_N =\{f_i, f_{ij}: i,j =1, \dots N\}$
is composed of functions satisfying \eqref{eq:condf1}-\eqref{eq:condf3}.
Moreover, we assume that associated to $\mathcal F = (\mathcal F_N)_{N \in \mathbb N}$ and to 
\begin{equation}\label{eq:dyn6}
\dot x_i(t)=f_i(\D x(t))+\sum_{j=1}^N f_{ij}(\D x(t))x_j(t),
\end{equation}
we can define a mean-field equation
\begin{equation}
  \frac{\partial \mu}{\partial t} + \nabla \cdot ( \mathcal H[\mathcal F](\mu) \mu) =0, \quad \mu(0)=\mu_0 \in \mathcal P_c(\mathbb R^d),
\end{equation}
such that for any compactly supported measure valued weak solutions
$\mu(t),\nu(t) \in C([0,T],\mathcal P_c(\mathbb R^d))$ of \eqref{genkin}
we have the following stability
\begin{equation}\label{stability2}
W_1(\mu(t), \nu(t)) \leq C(t) W_1(\mu(0), \nu(0)), \quad t \in [0,T],
\end{equation}
where $C(t)$ is a positive increasing function of $t$, independent of the dimension $d$. We further require that corresponding assumptions, including stability, hold for the projected system \eqref{eq:eul3} and kinetic equation \eqref{genkin3}.
Then we have the following approximation result:

\begin{theorem}\label{Thm:Kin}
Let us assume that $\mu_0 \in \mathcal P_c(\mathbb R^d)$ and there exist points
$\{x_1^0, \dots, x_N^0\} \subset \mathbb R^d$, for which the atomic measure
$\mu_0^N= \frac{1}{N} \sum_{i=1}^N \delta_{x_i^0}$ approximates
$\mu_0$ up to $\varepsilon>0$ in Monge-Kantorovich-Rubinstein distance, in the following sense
\begin{equation}\label{epskin}
W_1(\mu_0, \mu_0^N) \leq \varepsilon,\quad  N=\mathcal N^{\bar k(\varepsilon)}\mbox{ for } \bar k(\varepsilon) \leq d \mbox{ and } \bar k(\varepsilon) \to d \mbox{ for } \varepsilon \to 0.
\end{equation}
Requirement \eqref{epskin} is in fact called the \emph{delayed curse of dimension} as explained below in detail in Section \ref{delay}.
Depending on $\varepsilon>0$ we fix also
$$
k = k(\varepsilon) =\mathcal O(\varepsilon^{-2} \log(N)) = \mathcal O(\varepsilon^{-2} \log(\mathcal N) \bar k(\varepsilon)).
$$
Moreover, let $M:\mathbb R^d \to \mathbb R^k$ be a linear mapping
which is a \emph{continuous Johnson-Lindenstrauss embedding} as in \eqref{Cont-JL}
for continuous in time trajectories $x_i(t)$ of \eqref{eq:dyn6} with initial datum $x_i(0)=x_i^0$.
Let $\nu \in C([0,T], \mathcal P_c(\mathbb R^k))$ be the weak solution of
\begin{eqnarray}\label{genkin4}
&&\frac{\partial \nu}{\partial t} + \nabla \cdot ( \mathcal H[\mathcal F_M](\nu) \nu) =0, \\
&& \nu(0)= (\mu_0)_M \in \mathcal P_c(\mathbb R^k),\label{genkin5}
\end{eqnarray}
where $(\mu_0)_M =  M \# \mu_0$. Then
\begin{equation}\label{eq:approxkin}
W_1(\mu_M(t),\nu(t))\leq \mathcal C(t) \|M\| \varepsilon, \quad t \in [0,T],
\end{equation}
where $\mathcal C(t)$ is an increasing function of $t$, with $\mathcal C(0) >0$, which is at most polynomially growing with the dimension $d$.
\end{theorem}
\begin{proof}
Let us define $\nu^N(t)$ the solution to equation \eqref{genkin4} with initial datum $\nu^N(0)=(\mu_0^N)_M$, or, equivalently, thanks to \eqref{projatom}
$$
\nu^N(t) = \frac{1}{N} \sum_{i=1}^n \delta_{y_i(t)},
$$
where $y_i(t)$ is the solution of
\begin{eqnarray*}
    \dot y_i &=& f_i(\D' y) + \sum_{j=1}^N f_{ij}(\D' y) y_j \,,\qquad i=1,\dots, N\,, \nonumber\\
      y_i(0) &=& Mx_i^0 \,,\qquad i=1,\dots, N\,.
\end{eqnarray*}
We estimate 
$$
W_1(\mu_M(t), \nu(t)) \leq W_1(\mu_M(t),(\mu^N(t))_M) + W_1((\mu^N(t))_M, \nu^N(t))+ W_1(\nu^N(t), \nu(t)).
$$
By using the definition of push-forward \eqref{pf} and \eqref{epskin}, the first term can be estimated by
\begin{eqnarray*}
W_1(\mu_M(t),(\mu^N(t))_M) &=& \sup\{\langle \mu_M(t) - (\mu^N(t))_M,\varphi \rangle: {\rm Lip}(\varphi) \leq 1\} \\
&=& \sup\{\langle \mu(t) - \mu^N(t),\varphi(M \cdot) \rangle: {\rm Lip}(\varphi) \leq 1\}\\
&\leq& \|M\| W_1(\mu(t), \mu^N(t)) \leq \|M\| C(t) \varepsilon.
\end{eqnarray*}
We estimate now the second term
\begin{eqnarray*}
W_1((\mu^N(t))_M, \nu^N(t)) &=& \sup\{\langle (\mu^N(t))_M- \nu^N(t),\varphi \rangle: {\rm Lip}(\varphi) \leq 1\} \\
&=& \sup\{\frac{1}{N} \sum_{i=1}^N (\varphi(M x_i(t)) - \varphi(y_i(t))): {\rm Lip}(\varphi) \leq 1\}\\
&\leq& \frac{1}{N} \sum_{i=1}^N \| Mx_i(t) - y_i(t)\|_{\ell_2^k}. 
\end{eqnarray*}
We recall the uniform approximation of Theorem~\ref{Thm:Cont},
$$
\| Mx_i(t) - y_i(t)\|_{\ell_2^k} \leq D(t) \varepsilon\,,\qquad i=1,\dots, N,
$$
where $D(t)$ is the time-dependent function on the right-hand-side of \eqref{eq:unifstab2}. Hence
$$
W_1(\mu_M(t),(\mu^N(t))_M) \leq D(t) \varepsilon.
$$
We address now the upper estimate of the third term, by the assumed stability of the lower dimensional equation \eqref{genkin3}
\begin{eqnarray*}
W_1(\nu^N(t), \nu(t)) &\leq& C(t) W_1(\nu^N(0), \nu(0)) \\
&= & C(t)  W_1((\mu^N_0)_M, (\mu_0)_M) \\
&\leq& C(t) \|M\| W(\mu^N_0,\mu_0) \leq  C(t) \|M\| \varepsilon.
\end{eqnarray*}
We can fix $\mathcal C(t) = 2  C(t) \|M\| + D(t)$, and, as observed in Theorem \ref{thmcontJL}, we can assume without loss of generality that $\|M\|  \leq \sqrt{\frac{d}{k}}$. Hence, $\mathcal C(t)$ depends at most polynomially with respect to the dimension $d$.
\end{proof}

\subsection{Approximation of probability measures by atomic measures and optimal integration}\label{sec:appr}

In view of the fundamental requirement \eqref{epskin} in Theorem \ref{Thm:Kin},
given $\mu_0 \in \mathcal P_c(\mathbb R^d)$, we are interested to establish
an upper bound to the best possible approximation in Monge-Kantorovich-Rubinstein distance
by means of atomic measures $\mu^N_0 = \frac{1}{N}\sum_{i=0}^{N-1} \delta_{x_i^0}$
with $N$ atoms, i.e., 
\begin{eqnarray}
\mathcal E_N(\mu_0)&:=&\inf_{\mu^N_0= \frac{1}{N}\sum_{i=0}^{N-1} \delta_{x_i^0}} W_1(\mu_0,\mu^N_0) \label{optint} \\
&=& \inf_{\{x_0^0,\dots, x_{N-1}^0\} \subset \mathbb R^d} \sup  \big \{ | \int_{\mathbb R^d} \xi(x) d\mu_0(x) - \frac{1}{N} \sum_{i=0}^{N-1} \xi(x_i^0) |: \xi \in {\rm Lip}(\mathbb R^d),  {\rm Lip}(\xi) \leq 1 \big \}   \nonumber.
\end{eqnarray}
In fact, once we identify the optimal points $\{x_0^0,\dots, x_{N-1}^0\}$,
we can use them as initial conditions $x_i(0)=x_i^0$ for the dynamical system \eqref{eq:dyn6},
and by using the stability relationship \eqref{stability}, we obtain
\begin{equation}\label{stability_red}
W_1(\mu(t), \mu^N(t)) \leq C(T) W_1(\mu_0, \mu^N_0), \quad t \in [0,T] \,,
\end{equation}
where $\mu^N(t) = \frac{1}{N} \sum_{i=0}^{N-1} \delta_{x_i(t)}$, 
meaning that the solution of the partial differential equation \eqref{genkin}
keeps optimally close to the particle solution of \eqref{eq:dyn6} also for successive time $t>0$.
Note that estimating \eqref{optint} as a function of $N$
is in fact a very classical problem in numerical analysis well-known as {\it optimal integration}
with its high-dimensional behaviour being a relevant subject of the field of {\it Information Based Complexity} \cite{NW,TWW}.

The numerical integration of Lipschitz functions with respect to the Lebesgue measure and the study
of its high-dimensional behaviour goes back to Bakhvalov \cite{Bakh}, but much more is known nowadays.
We refer to \cite{GL} and \cite{Gruber} for the state of the art of quantization of probability distributions.

The scope of this section is to recall some facets of these estimates and to reformulate them in terms of $W_1$ and ${\mathcal E}_N$. 
We emphasize that here and in what follows, we consider generic compactly supported probability measures $\mu$, not necessarily absolutely continuous with respect to the Lebesgue measure. 
We start first by assuming $d=1$, i.e., we work with a univariate measure $\mu \in \mathcal P_c(\mathbb R)$ with 
support $\supp \mu \subset [a, b]$ and $\sigma :=  b - a >0$. We define the points $x_0,\dots,x_{N-1}$ 
as the {\it quantiles} of the probability measure $\mu$, i.e., $x_0:=a$ and
\begin{equation}\label{quantiles}
\frac{i}{N} = \int_{-\infty}^{x_i} d\mu(x), \quad i=1,\dots, N-1.
\end{equation}
This is notationally complemented by putting $x_N:=b$.
Note that by definition $\int_{x_{i}}^{x_{i+1}} d\mu(x) =\frac{1}{N}, i=0,\dots, N-1$, and we have
\begin{eqnarray}
\notag\left|\int_{\mathbb R} \xi(x) d\mu(x) - \frac{1}{N} \sum_{i=0}^{N-1} \xi(x_i) \right | &=& \left | \sum_{i=0}^{N-1} \int_{x_i}^{x_{i+1}} (\xi(x)- \xi(x_i)) d\mu(x) \right |\\
\label{quantiles2}&\leq& \sum_{i=0}^{N-1} \int_{x_i}^{x_{i+1}}\left |\xi(x)- \xi(x_i)\right | d\mu(x)  \\
\notag&\leq& \frac{{\rm Lip}(\xi)}{N} \sum_{i=0}^{N-1} (x_{i+1}-x_i) = \frac{{\sigma} \rm  Lip(\xi)}{N}.
\end{eqnarray}
Hence it is immediate to see that 
$$
\mathcal E_N(\mu)=\inf_{\mu^N= \frac{1}{N}\sum_{i=0}^{N-1} \delta_{x_i^0}} W_1(\mu,\mu^N) \leq \frac{{\sigma}}{N}.
$$
We would like to extend this estimate to higher dimension $d >1$. However, for multivariate measures $\mu$ there is no such an
easy upper bound, see  \cite{GL} and \cite{Gruber} for very general statements, and for the sake of simplicity we restrict here the class of measures $\mu$ to certain special cases.
As a typical situation, we address tensor product measures and sums of tensor products. 
\begin{lemma}\label{lem:W:tensors}
Let $\mu^1,\dots,\mu^d\in \mathcal P_1(\R)$ with $W_1(\mu^j,\mu^{j,N_j})\le \varepsilon_j,$ $j=1,\dots, d$
for some $N_1,\dots,N_d\in\N$, $\varepsilon_1,\dots,\varepsilon_d >0$ and $\mu^{j,N_j}:=\frac{1}{N_j}\sum_{i=0}^{N_j-1}\delta_{x^j_{i}}$.
Let $N=\prod_{i=1}^d N_i$. Then
\begin{equation*}
W_1(\mu^1\otimes\dots\otimes \mu^d,\mu^{N})\le \sum_{j=1}^d \varepsilon_j,
\end{equation*}
where
\begin{equation*}
\mu^N:=\frac{1}{N}\sum_{x\in X}\delta_x\quad\text{and}\quad X:=\prod_{j=1}^d\{x^j_0,\dots,x^j_{N_j-1}\}.
\end{equation*}
\end{lemma}
\begin{proof}
The proof is based on a simple argument using a telescopic sum. For $j=1,\dots,d+1$ we put
\begin{align*}
V_j&:=\frac{1}{\prod_{i=j}^d N_{i}}\sum_{i_{j}=0}^{N_{j}-1}\dots\sum_{i_d=0}^{N_d-1}\int_{\R^{j-1}}
\xi(x_1,\dots,x_{j-1},x^{j}_{i_{j}},\dots,x^d_{j_d}) d\mu^1(x_1)\dots d\mu^{j-1}(x_{j-1}).
\end{align*}
Of course, if $j=1$, then the integration over $\R^{j-1}$ is missing and if $j=d+1$ then the summation becomes empty.
Now
$$
\int_{\R^d}\xi(x) d\mu(x)-\frac{1}{\prod_{i=1}^d N_i}\sum_{i_1=0}^{N_1-1}\dots\sum_{i_d=0}^{N_d-1}\xi(x^{1}_{i_{1}},\dots,x^d_{i_d})
=\sum_{j=1}^{d}(V_{j+1}-V_j)
$$
together with the estimate $|V_{j+1}-V_j|\le \varepsilon_j$ finishes the proof.
\end{proof}

Lemma \ref{lem:W:tensors} says, roughly speaking, that the tensor products of sampling points of univariate
measures are good sampling points for the tensor product of the univariate measures.
Next lemma deals with sums of measures.

\begin{lemma}\label{lem:W:sums} Let $\mu_1,\dots,\mu_L\in \mathcal{P}_1(\R^d)$ with 
$W_1(\mu_l,\mu_l^{N})\le \varepsilon_l$, $l=1,\dots, L$
for some $N\in\N$, $\varepsilon_1,\dots,\varepsilon_L >0$ and $\mu_l^{N}:=\frac{1}{N}\sum_{i=0}^{N-1}\delta_{x_{l,i}}$.
Then
\begin{equation*}
W_1\Bigl(\frac{\mu_1+\dots+ \mu_L}{L},\mu^{LN}\Bigr)\le \frac{1}{L}\sum_{l=1}^L \varepsilon_l,
\end{equation*}
where
\begin{equation*}
\mu^{LN}:=\frac{1}{LN}\sum_{x\in X}\delta_x=\frac{1}{L}\sum_{l=1}^L \mu_l^{N}\quad\text{and}\quad 
X:=\bigcup_{l=1}^L\{x_{l,0},\dots,x_{l,N-1}\}.
\end{equation*}
\end{lemma}
\begin{proof}
We use the homogeneity of the Monge-Kantorovich-Rubinstein distance $W_1(a\mu,a\nu)=aW_1(\mu,\nu)$ for $\mu,\nu\in\mathcal P_1(\R^d)$
and $a\ge 0$ combined with its subadditivity 
$W_1(\mu_1+\mu_2,\nu_1+\nu_2)\le W_1(\mu_1,\nu_1)+W_1(\mu_2,\nu_2)$ for $\mu_1,\mu_2,\nu_1,\nu_2\in\mathcal P_1(\R^d)$. We obtain
\begin{equation*}
W_1\Bigl(\frac{\mu_1+\dots+ \mu_L}{L},\frac{\mu_1^{N}+\dots+\mu_L^{N}}{L}\Bigr)\le \frac{1}{L}\sum_{l=1}^L W_1(\mu_l,\mu_l^{N})\le \frac{1}{L}\sum_{l=1}^L \varepsilon_l.
\end{equation*}
\end{proof}

Next corollary follows directly from Lemma \ref{lem:W:tensors} and Lemma \ref{lem:W:sums}.
\begin{corollary}\label{cor1} (i) Let $\mu^1,\dots,\mu^d\in \mathcal P_1(\R)$ and $N_1,\dots,N_d\in\N$. 
Then 
\begin{equation*}
\mathcal E_{N}(\mu^1\otimes\dots\otimes\mu^d)\le \sum_{j=1}^d \mathcal E_{N_j}(\mu^j),\quad \text{where} \quad N:=N_1\cdots N_d.
\end{equation*}
(ii) Let $\mu_1,\dots,\mu_L\in \mathcal{P}_1(\R^d)$ and $N\in\N$. Then
\begin{equation*}
{\mathcal E}_{LN}\Bigl(\frac{\mu_1+\dots+ \mu_L}{L}\Bigr)\le \frac{1}{L}\sum_{l=1}^L {\mathcal E}_N(\mu_l).
\end{equation*}
\end{corollary}

\subsection{Delayed curse of dimension}\label{delay}

Although Lemma \ref{lem:W:tensors}, Lemma \ref{lem:W:sums} and Corollary \ref{cor1} give some estimates of
the Monge-Kantorovich-Rubinstein distance between general and atomic measures, the number of atoms needed may still be too large to allow the assumption \eqref{epskin} in Theorem \ref{Thm:Kin} to be fulfilled.
Let us for example consider the case, where $\mu^1=\dots=\mu^d$ in Lemma \ref{lem:W:tensors} and $\varepsilon_1=\dots=\varepsilon_d=:\varepsilon.$
Then, of course, $N_1=\dots=N_d=:{\mathcal N}$ and we observe, that the construction given in Lemma \ref{lem:W:tensors}
gives an atomic measure, which approximates $\mu$ up to the error $d\varepsilon$ using ${\mathcal N}^d$ atoms,
hence with an exponential dependence on the dimension $d$. This effect is another instance of the well-known 
phenomenon of the {\it curse of dimension}.

However, in many real-life high-dimensional applications the objects of study (in our case the measure $\mu\in{\mathcal P}_c(\R^d)$)
concentrate along low-dimensional subspaces (or, more general, along low-dimensional manifolds)  \cite{BN01,BN03,CLLMNWZ05a,CLLMNWZ05b,CL06}. 
The number of atoms necessary to approximate these measures behaves in a much better way, allowing the application of  \eqref{epskin} and Theorem \ref{Thm:Kin}. To clarify this effect, let us consider
$\mu=\mu^1\otimes\dots\otimes \mu^d$ with $\supp\mu^j\subset [a_j,b_j]$ and define $\sigma_j=b_j-a_j$. Let us assume, that
$\sigma_1\ge \sigma_2\ge \dots\ge \sigma_d > 0$ is a rapidly decreasing sequence. Furthermore, let $\varepsilon>0$.
Then we define $\bar k:=\bar k(\varepsilon)$ to be the smallest natural number, such that
$$
\sum_{k=\bar k(\varepsilon)+1}^{d}\sigma_k\le \varepsilon/2
$$
and put $N_k=1$ for $k\in\{\bar k(\varepsilon)+1,\dots,d\}$. The numbers $N_1=\dots=N_{\bar k(\varepsilon)}={\mathcal N}$ are chosen large enough so that
$$
\frac{1}{\mathcal N}\sum_{k=1}^{\bar k(\varepsilon)} {\sigma_k}\le \varepsilon/2.
$$
Then Lemma \ref{lem:W:tensors} together with \eqref{quantiles} state that there is an atomic measure $\mu^N$ 
with $N={\mathcal N}^{\bar k(\varepsilon)}$ atoms, such that
\begin{equation}\label{delayed}
W_1(\mu,\mu^N)\le \sum_{k=1}^d \frac{\sigma_k}{N_k}\le \varepsilon/2+\varepsilon/2.
\end{equation}
Hence, at the cost of assuming that the tensor product measure $\mu$ is concentrated along a $\bar k(\varepsilon)$-dimensional coordinate subspace,
we can always approximate the measure $\mu$ with accuracy $\varepsilon$ by using an atomic measure supported on points whose number 
depends exponentially on $\bar k=\bar k(\varepsilon) \ll d$.
However, if we liked to have $\varepsilon \to 0$, then $\bar k(\varepsilon) \to d$ and again we are falling under the curse of dimension.
This delayed kicking in of the need of a large number of points for obtaining high accuracy in the approximation \eqref{delayed}
is in fact the so-called {\it delayed curse of dimension}, expressed by assumption \eqref{epskin}, a concept introduced first by Curbera in  \cite{C00},
in the context of optimal integration with respect to Gaussian measures in high dimension.

Let us only remark, that the discussion above may be easily extended (with help of Lemma \ref{lem:W:sums}) 
to sums of tensor product measures. In that case we obtain as atoms the so-called \emph{sparse grids}, cf. \cite{BG}.
Using suitable change of variables, one could also consider measures concentrated
around (smooth) low-dimensional manifolds, but this goes beyond the scope of this work, see \cite{GL} for a broader discussion.




\subsubsection*{Acknowledgments}
We acknowledge the financial support provided by the START award
``Sparse Approximation and Optimization in High Dimensions'' no. FWF~Y~432-N15
of the Fonds zur F\"orderung der wissenschaftlichen Forschung (Austrian Science Foundation).
We would also like to thank Erich Novak for a discussion about optimal integration
and for pointing us to some of the references given in Section \ref{sec:appr}.

\end{document}